\documentclass{conm-p-l}

\usepackage{epic,eepic}
\usepackage{verbatim}

\newtheorem{thm}[equation]{Theorem}
\newtheorem{cor}[equation]{Corollary}
\newtheorem{lem}[equation]{Lemma}
\newtheorem{prop}[equation]{Proposition}

\theoremstyle{definition}
\newtheorem{defn}[equation]{Definition}

\theoremstyle{remark}
\newtheorem{rem}[equation]{Remark}
\newtheorem{ex}[equation]{Example}
\newtheorem{notation}{Notation}

\numberwithin{equation}{section}

\newcommand{\thmref}[1]{Theorem~\ref{#1}}
\newcommand{\propref}[1]{Proposition~\ref{#1}}
\newcommand{\lemref}[1]{Lemma~\ref{#1}}
\newcommand{\corref}[1]{Corollary~\ref{#1}}
\newcommand{\defref}[1]{Definition~\ref{#1}}
\newcommand{\remref}[1]{Remark~\ref{#1}}
\newcommand{\exref}[1]{Example~\ref{#1}}

\newcommand{\secref}[1]{Section~\ref{#1}}

\newcommand\1{{\text{\bf 1}}}
\renewcommand\a{\alpha}
\newcommand\BB{{\mathcal B}}
\renewcommand\b{\beta}

\newcommand\card{\operatorname{card}}
\newcommand\CC{{\mathcal C}}
\newcommand\D{\Delta}
\renewcommand\d{\delta}
\newcommand\diag{\operatorname{diag}}
\newcommand\E{{\mathbf{E}}}
\newcommand\e{\varepsilon}

\newcommand\F{\Phi}
\newcommand\FF{\mathcal{F}}
\newcommand\f{\varphi}
\newcommand\ft[1]{{\parbox{100truemm}{#1}}}

\newcommand\g{\gamma}

\renewcommand\H{{\mathbf{H}}}
\newcommand\HH{{\mathcal{H}}}

\newcommand\Iso{\operatorname{Iso}}
\renewcommand\l{\lambda}

\newcommand\M{{\mathfrak M}}
\newcommand\m{{\mathfrak m}}

\newcommand\Mu{{\mathrm M}}

\newcommand\OO{\mathcal O}
\newcommand\om{\omega}
\newcommand\ov{\overline}

\renewcommand\Pr{\mathcal{P}}
\renewcommand\part{\partial}
\newcommand\R{\mathbb R}
\newcommand\RR{\mathcal R}
\newcommand\s{\sigma}
\newcommand\supp{\operatorname{supp}}
\newcommand\TT{\mathcal{T}}
\renewcommand\th{\theta}

\renewcommand\th{\theta}
\newcommand\toto{\mathop{\;\longrightarrow\;}}
\newcommand\un{\underline}
\newcommand\wh{\widehat}
\newcommand\wt{\widetilde}

\newcommand\XX{\mathcal X}
\newcommand\x{{\boldsymbol x}}

\newcommand\Z{{\mathbb Z}}

\begin{document}

\title{Boundary amenability of hyperbolic spaces}

\author{Vadim A. Kaimanovich}

\address{CNRS UMR 6625, IRMAR, Universit\'{e} Rennes-1, Campus Beaulieu, 35042
Rennes, France}

\email{kaimanov@univ-rennes1.fr}

\thanks{
The paper was mostly written during author's stay at the Max-Planck-Institut
(Bonn) in 2002. Its financial support and the excellent working conditions are
gratefully acknowledged. I would also like to thank the organizers of the First
JAMS Symposium ``Discrete Geometric Analysis'', where these results were
presented, and, in particular, Professor Motoko Kotani, for creating an ideal
atmosphere at this very fruitful meeting. I am grateful to Theo B\"{u}hler and
to the anonymous referee for a number of useful comments and suggestions.}



\begin{abstract}
It is well-known that a Kleinian group is amenable if and only if it is
elementary. We establish an analogous property for equivalence relations and
foliations with Gromov hyperbolic leaves: they are amenable if and only if they
are elementary in the sense that one can assign (in a measurable way) to any
leaf a finite subset of its hyperbolic boundary (as in the group case, such
subsets cannot actually contain more than 2 points). The analogous result for
actions of word hyperbolic groups with a quasi-invariant measure is that such
an action is amenable if and only if it factorizes through the hyperbolic
boundary or its symmetric square.

A byproduct of our approach is a proof of boundary amenability for isometry
groups of Gromov hyperbolic spaces under the assumption that either the space
is exponentially bounded or it is CAT($-1$) and the group has a finite critical
exponent. We also give examples showing that without these assumptions the
boundary amenability may fail.
\end{abstract}

\maketitle

\section*{Introduction}

The class of \emph{amenable groups} is, from the analytical point of view, the
most natural generalization of the class of compact groups. Amenable groups are
those which admit an \emph{invariant mean} (rather than an invariant
probability measure, which is the case for compact groups). There are many
different equivalent definitions of amenability. Probably, the most
constructive definition (used for verifying amenability of a given group) is
the one formulated in terms of existence of \emph{approximatively invariant
sequences of probability measures} on the group (\emph{Reiter's condition}),
whereas one of the main applications of amenability is the \emph{fixed point
property} for affine actions of amenable groups on compact spaces.

It turns out that non-amenable groups may still have actions which look like
actions of amenable groups. This observation led Zimmer \cite{Zimmer77},
\cite{Zimmer78} to introduce the notion of an \emph{amenable action}. In the
same way as with groups, there are several definitions of an amenable action.
In particular, amenable actions can be characterized both in terms of a fixed
point property (this was the original definition of Zimmer) and in terms of
existence of a sequence of approximatively equivariant maps from the action
space to the space of probability measures on the group (this is an analogue of
Reiter's condition). Yet another generalization is the notion of amenability
for \emph{equivalence relations} and \emph{foliations}. Actually, all these
objects can be considered as \emph{measured groupoids}, and the notion of
amenability in each particular case is a specialization of a general notion of
an \emph{amenable groupoid}, see \cite{Anantharaman-Renault00},
\cite{Corlette-Lamoneda-Iozzi02}.

\medskip

It is well-known that a Kleinian group is \emph{amenable} if and only if it is
\emph{elementary}. We establish an analogous property for \emph{equivalence
relations and foliations with Gromov hyperbolic leaves}: under a suitable
bounded geometry condition they are amenable if and only if they are elementary
in the sense that one can assign (in a measurable way) to any leaf a finite
subset of its hyperbolic boundary (actually, such subsets cannot consist of
more than 2 points). As particular cases it contains and generalizes earlier
results on \emph{amenability of boundary actions} of isometry groups of
hyperbolic spaces and on \emph{amenability of the stable foliation} of the
geodesic flow.

The main geometric ingredients of our approach --- (i) \emph{boundary
convergence of geodesics} in hyperbolic spaces (in various guises), and (ii)
the notion of the \emph{barycenter} for hyperbolic spaces --- are well-known
and have been used before (e.g., see \cite{Adams96}, \cite{Burger-Mozes96}).
However, working with hyperbolic equivalence relations requires developing a
new technique (based on the notion of a \emph{measurable bundle of Banach
spaces}) which allows one to deal with \emph{measurable systems of} (possibly
different) \emph{leafwise hyperbolic boundaries}. Another technical novelty in
this paper is the use of simple explicit constructions of approximatively
invariant sequences of measures as \emph{Cesaro averages for sandwiched
sequences of sets with subexponential growth} (\thmref{th:lambda} and
\lemref{lem:sandwich}) and \emph{Cesaro averages of pre-Patterson measures
along geodesics} (\thmref{th:lambda-CAT}).

\medskip

The paper has the following structure. In \secref{sec:hyperbolic} we collect
the geometric properties of Gromov hyperbolic spaces used in the sequel. These
are \propref{prop:sausage} on boundary convergence of geodesic rays,
\propref{pr:allbary} on existence of (quasi)-barycenter sets and Theorems
\ref{th:lambda}, \ref{th:lambda-CAT} on existence of approximatively invariant
boundary maps.

The main technical tool used to ensure measurability of leafwise application of
these constructions in the setup of graphed equivalence relations with
hyperbolic leaves is the measurable Banach bundle of leafwise functions with a
continuous boundary extension (\secref{sec:equiv relations}). By using this
bundle we define measurable systems of boundary measures and show that a
hyperbolic equivalence relation admits an invariant boundary system of
probability measures if and only if it is elementary (\thmref{th:3types}).
Further we establish an analogous result for measurable foliations
(\thmref{thm:classes of foliations}) by reducing them to equivalence relations.

We begin \secref{sec:amenability} with illustrating the key geometric
ingredients of the relationship between amenability and hyperbolicity, namely,
boundary convergence of geodesics and existence of quasi-barycenter sets, on
the example of coincidence of the class of amenable groups of isometries of
Gromov hyperbolic spaces with the class of elementary groups (under suitable
bounded geometry type conditions). Further we proceed to establish analogues of
this result for equivalence relations, foliations and actions
(\thmref{th:equivmain}, \thmref{th:foliations main} and \thmref{th:actions
mains}, respectively) by taking stock of the results from the previous sections
and following precisely the same scheme: amenability by the fixed point
property implies existence of an invariant system of boundary measures, and
therefore elementarity, whereas elementarity leads to existence of
approximatively invariant sequences of probability measures, i.e., to
amenability. As a byproduct we give a proof of the topological amenability of
the boundary action of the group of isometries of a hyperbolic space under the
assumption that either the space is exponentially bounded or the space is
CAT($-1$) and the group has a finite critical exponent (\thmref{th:univ} and
\thmref{th:univ-CAT}). Finally, we give examples showing that without these
assumptions the boundary amenability (even in the weakest form) may well fail.

\section{Asymptotic geometry of hyperbolic spaces} \label{sec:hyperbolic}

\subsection{Hyperbolic spaces}

A detailed discussion of the notion of $\d$-hyperbo\-li\-ci\-ty and of the
associated structures can be found in the seminal work of Gromov
\cite{Gromov87} and in the notes \cite{Ghys-delaHarpe90}. Below are listed some
of the definitions and properties used later on.

We shall choose the definition based on the \emph{Rips condition}: a
non-compact complete proper geodesic metric space $\XX$ is
\emph{$\d$-hyperbolic} (with $\d\ge 0$) if each of the sides of any geodesic
triangle is contained in the $\d$-neighbourhood of the union of the other two
sides (see, for instance, \cite[Proposition 2.21]{Ghys-delaHarpe90} for a list
of other equivalent definitions). The minimal number $\d$ with this property is
the \emph{hyperbolicity constant} of $\XX$. A graph is called $\d$-hyperbolic
if the associated 1-complex with length 1 edges is $\d$-hyperbolic. Usually, we
shall not be concerned with the precise value of the hyperbolicity constant,
and call the above spaces just \emph{hyperbolic}.

\medskip

For the rest of this Section we shall fix a hyperbolic space $\XX$ with metric
$d$ and the hyperbolicity constant $\d$. Denote by
$$
(y|z)_x = \frac12 \bigl[ d(x,y) + d(x,z) - d(y,z) \bigr] \;, \qquad x,y,z\in\XX
$$
the \emph{Gromov product} on $\XX$. Then by \cite[Lemme 2.17 and Proposition
2.21]{Ghys-delaHarpe90} for any geodesic segment $[y,z]$ joining the points $y$
and $z$
\begin{equation} \label{eq:4d}
d(x,[y,z]) - 4\d \le (y|z)_x \le d(x,[y,z]) \;.
\end{equation}

The functions on $\XX\times\XX$
\begin{equation} \label{eq:varrho}
\varrho_x(y,z) =
  \begin{cases}
    0 \;, & y=z \;, \\
    e^{-(y|z)_x} \;, & \text{otherwise\;,}
  \end{cases}
\end{equation}
determine a uniform structure ``at infinity''on the space $\XX$ (which does not
depend on the choice of the point $x$). The completion $\ov\XX$ of $\XX$ with
respect to this uniform structure is called the \emph{hyperbolic
compactification} of $\XX$. Its boundary is denoted by $\part
\XX=\ov\XX\setminus\XX$. The action of the group of isometries $\Iso(\XX)$
extends from $\XX$ to a continuous action on $\part\XX$.

Any geodesic ray $\xi$ converges in the hyperbolic compactification to a
boundary point $\xi(\infty)\in\part\XX$. The definition of $\d$-hyperbolicity
is easily seen to imply the following property (which is, for instance,
essentially contained in \cite[Proposition~7.2]{Ghys-delaHarpe90}):

\begin{prop} \label{prop:sausage}
There exists a constant $C_1=C_1(\d)$ with the property that for any
$\d$-hyperbolic space $\XX$ and any two geodesic rays $\xi_1,\xi_2$ in $\XX$
with $\xi_1(\infty)=\xi_2(\infty)$ there is a number $T=T(\xi_1,\xi_2)$ such
that $|T|\le d(\xi_1(0),\xi_2(0))$ and
$$
d(\xi_1(t),\xi_2(t+T))\le C_1  \qquad \forall\, t\ge d(\xi_1(0),\xi_2(0)) \;.
$$
\end{prop}

Although the functions $\varrho_x$ \eqref{eq:varrho} do not, in general,
satisfy the triangle inequality, they are ``almost'' metrics. Namely, for
$a=\frac1{15\d}$ the inner metric
\begin{equation} \label{eq:rho}
 \rho_x(y,z)
 = \inf \left\{ \sum_{i=1}^n [\varrho_x(x_{i-1},x_i)]^a :
 n \ge 1\,,\; y=x_0, x_1, \dots, x_n=z \in \XX \right\}
\end{equation}
determined by the function $\varrho_x^a$ satisfies the inequality
\begin{equation} \label{eq:1/2}
\frac12(\varrho_x)^a \le \rho_x \le (\varrho_x)^a
\end{equation}
(this is an adaptation of \cite[Proposition~22.8]{Woess00}). So, the hyperbolic
compactification $\ov\XX$ is the completion of $\XX$ with respect to any of the
metrics $\rho_x$.

\begin{notation}
Below we shall use the notation $B(o,r)$ (resp., $S(o,r)$) for the $r$-ball
(resp., the $r$-sphere) in $\XX$ of the hyperbolic metric $d$ around a point
$o\in\XX$, and $B_x(o,r)$ (resp., $S_x(o,r)$) for the $r$-ball (resp., the
$r$-sphere) in $\ov\XX$ of the metric $\rho_x$ around a point $o\in\ov\XX$. We
denote the space of probability measures on a topological space $X$ by
$\Pr(X)$, and the normalized restriction of a measure $m$ to a positive measure
set $A$ by $m_A$.
\end{notation}

\subsection{The barycenter set of a boundary measure}

Denote by
\begin{equation} \label{eq:bz}
 \b_z(x,y)
 = d(y,z) - d(x,z)
 = d(x,y) - 2(y|z)_x \qquad x,y\in\XX
\end{equation}
the \emph{distance cocycle} associated with a point $z\in\XX$. The triangle
inequality for the metric $d$ and formula \eqref{eq:1/2} imply

\begin{prop} \label{pr:estimates bz}
For any $x,y,z\in\XX$ the distance cocycle satisfies the inequalities
$$
|\b_z(x,y)| \le d(x,y) \;,
$$
and $($assuming $y\neq z)$
$$
d(x,y) + \frac2a\log\rho_x(y,z)
 \le  \b_z(x,y)
 \le d(x,y) + \frac2a\log\rho_x(y,z) + C_2\;,
$$
where
$$
C_2 = C_2(\d) = \frac{2\log2}a = (30\log 2) \d \;.
$$
\end{prop}

\begin{cor} \label{cor:limsupinf}
For any $\g\in\part\XX$ and $x,y\in\XX$
$$
\limsup_{z\to\g} \b_z(x,y) - \liminf_{z\to\g} \b_z(x,y) \le C_2 \;.
$$
\end{cor}

Put
$$
 \ov\b_\g(x,y) = \limsup_{z\to\g} \b_z(x,y)
 \qquad x,y\in\XX,\,\g\in\part\XX \;.
$$
As it follows from \eqref{eq:bz}, the functions $\ov\b_\g$ are Lipschitz with
respect to each of the arguments. Although $\ov\b_\g$ are not, generally
speaking, cocycles, they still satisfy the cocycle identity with a uniformly
bounded error (i.e., they are \emph{quasi-cocycles}). More precisely,
\propref{pr:estimates bz} and \corref{cor:limsupinf} imply

\begin{prop} \label{pr:quasicocycle}
For any $\g\in\part\XX$ the function $\ov\b_\g$ has the following properties:
\begin{itemize}
\item[(i)]
$\ov\b_\g$ is ``jointly Lipschitz'', i.e., $|\ov\b_\g(x,y)|\le d(x,y)$ for all
$x,y\in\XX$, in particular, $\ov\b_\g(x,x)\equiv 0\;;$
\item[(ii)]
$0\le\ov\b_\g(x,y) + \ov\b_\g(y,x) \le C_2$ for all $x,y\in\XX\;;$
\item[(iii)]
$0 \le  \ov\b_\g(x,y) + \ov\b_\g(y,z) + \ov\b_\g(z,x) \le 2C_2$ for all
$x,y,z\in\XX\;.$
\end{itemize}
\end{prop}

\begin{rem}
Our definition of $\ov\b_\g$ is different from the definition of the ``Busemann
function'' in \cite[Chapitre 8]{Ghys-delaHarpe90} (where the point $z$ was
allowed to converge to $\g$ along geodesic rays only).
\end{rem}

For any probability measure $\l\in\Pr(\part\XX)$ put
$$
\BB_\l(x,y) = \int \ov\b_\g(x,y) \, d\l(\g) \;.
$$
The Lipschitz property of the functions $\b_\g$ implies that $\BB_\l$ are also
Lipschitz with respect to each of the arguments. Moreover, as it follows from
\propref{pr:quasicocycle}, the functions $\BB_\l$ are jointly Lipschitz, and
they are quasi-cocycles with the same constants as in
\propref{pr:quasicocycle}. In particular,
\begin{equation} \label{eq:Bxx}
\bigl| \BB_\l(x,y) - \BB_\l(x',y) - \BB_\l(x,x') \bigr| \le 3 C_2 \qquad
\forall\,x,x',y\in\XX,\,\l\in\Pr(\part\XX) \;.
\end{equation}

\begin{prop} \label{pr:b-infinity}
If a sequence of probability measures $\l_n\in\Pr(\XX)$ is convergent in the
weak$^*$ topology to a measure $\l\in\Pr(\XX)$, then
$$
\left| \BB_\l(x,y) - \limsup_{n\to\infty} \int \b_z(x,y)\,d\l_n(z) \right| \le
C_2 \qquad \forall\,x,y\in\XX \;.
$$
\end{prop}

\begin{proof}
As it follows from \propref{pr:estimates bz}, if $x\neq y$ then both terms in
the left-hand side belong to the interval
$$
\left[ d(x,y) + \frac2a\int\log\rho_x(y,\g)\,d\l(\g), \;  d(x,y) +
\frac2a\int\log\rho_x(y,\g)\,d\l(\g) +C_2 \right] \;.
$$
\end{proof}

\begin{prop} \label{pr:escape}
If a measure $\l\in\Pr(\part\XX)$ does not have atoms of weight at least
$\displaystyle\frac12$, then for any $x\in X$
$$
\lim_y \BB_\l(x,y) = \infty
$$
if $y\to\infty$ in the space $\XX$ $($i.e., $d(x,y)\to\infty)$.
\end{prop}

\begin{proof}
Without loss of generality we may assume that $y\to\g_0\in\part\XX$. Choose
$\e>0$ in such a way that  $\l B_x(\g_0,\e) < \displaystyle\frac12$, and, by
using the inequalities from \propref{pr:estimates bz} (which are satisfied for
the functions $\ov\b_\g$ as well), estimate $\ov\b_\g(x,y)$ from below as
$$
\ov\b_\g(x,y)\ge
  \begin{cases}
    -d(x,y) \;, & \g\in B_x(\g_0,\e) \;, \\
     d(x,y) + \frac2a\log\rho_x(y,\g) \;, &
     \g\in\part\XX\setminus B_x(\g_0,\e) \;.
  \end{cases}
$$
Since $y\to\g_0$, we may assume that $\rho_x(y,\g)\ge\e/2$ for any
$\g\in\part\XX\setminus B_x(\g_0,\e)$, whence integrating the above estimate
yields the claim.
\end{proof}

\begin{rem} \label{rem:Busemann}
If $\XX$ is a CAT($-1$) space, then all functions $\ov\b_\g$ and $\BB_\l$ are
actually cocycles, so that under conditions of \propref{pr:escape} the function
$\BB_\l(x,\cdot)$ attains its minimum on a compact subset of $\XX$ (the
\emph{barycenter set} of $\l$) which does not depend on the choice of the point
$x\in\XX$, see \cite{Douady-Earle86}, \cite{Besson-Courtois-Gallot96},
\cite{Burger-Mozes96}. Due to the usual technical complications connected with
replacing the CAT($-1$) property with $\d$-hyperbolicity, this is no longer
true for a general $\d$-hyperbolic space, which is why we have to ``quasify''
the definition of the barycenter set.
\end{rem}

\begin{defn} \label{def:baryset}
Given a number $r>0$ the set
$$
\CC(\l,x,r) = \left\{y\in\XX: \BB_\l(x,y) \le r+ \inf_{z\in\XX} \BB_\l(x,z)
\right\}
$$
is called the \emph{quasi-barycenter set} of the measure $\l\in\Pr(\part\XX)$
with respect to the point $x\in\XX$.
\end{defn}

Inequality \eqref{eq:Bxx} and \propref{pr:escape} imply

\begin{prop} \label{pr:bary}
If a measure $\l\in\Pr(\part\XX)$ does not have atoms of weight at least
$\displaystyle\frac12$, then the infimum in \defref{def:baryset} is finite, the
quasi-barycenter set \linebreak $\CC(\l,x,r)$ is compact for any $x\in\XX$ and
$r>0$, and
$$
\CC(\l,x,r) \subset \CC(\l,x',r+6C_2) \qquad \forall\,x,x'\in\XX \;.
$$
\end{prop}

\begin{defn} \label{def:allbaryset}
We shall call the closure
$$
\CC(\l,r) = \ov{\bigcup_{x\in\XX} \CC(\l,x,r)} \;,
$$
the \emph{quasi-barycenter set} of the measure $\l\in\Pr(\part\XX)$.
\end{defn}

\propref{pr:bary} then implies (cf. \cite[Proposition 5.1]{Adams96}):

\begin{prop} \label{pr:allbary}
If a measure $\l\in\Pr(\part\XX)$ does not have atoms of weight at least
$\displaystyle\frac12$, then for any $r>0$ the quasi-barycenter set $\CC(\l,r)$
is compact, and the map $\l\mapsto\CC(\l,r)$ is $\Iso(\XX)$-equivariant.
\end{prop}

\subsection{Bounded geometry conditions}

Below we shall use several bounded geometry type conditions on a metric space
$\XX$.

\begin{defn} \label{def:tempered}
We say that a Radon measure $m$ on a metric space $\XX$ is
\emph{$(r_1,r_2)$-tempered} (with $0<r_1\le r_2\le\infty$) if for any $r\in
[r_1,r_2]$ the measures of closed $r$-balls in $\XX$ are uniformly bounded from
above and bounded away from $0$, i.e., there exist constants $0<\un C(r)\le \ov
C(r) <\infty,\;r\in[r_1,r_2]$ such that
$$
\un C(r) \le m B(x,r) \le \ov C(r) \qquad \forall\,x\in\XX\;.
$$
A measure is \emph{tempered} if it is $(r,\infty)$-tempered for a certain $r$.
A metric space $\XX$ is \emph{tempered} if it carries an $\Iso(\XX)$-invariant
tempered measure.
\end{defn}

\begin{defn}[cf. \cite{Adams96}] \label{def:at most}
A metric space $\XX$ is \emph{exponentially $\rho$-bounded} (with $\rho>0$) if
there exists a constant $a=a(\rho)$ such that, for every $x\in\XX$ and every
$r>0$, the ball $B(x,r)$ in $\XX$ can contain at most $a^r$ pairwise disjoint
balls of radius $\rho$. We shall say that $\XX$ is \emph{exponentially bounded}
if it is exponentially $\rho$-bounded for a certain $\rho$.
\end{defn}

\begin{rem}
Our terminology differs from that of \cite{Adams96}. What we call
``exponentially bounded'' is ``at most exponential'' in \cite{Adams96}.
\end{rem}

\begin{thm} \label{thm:measure}
A proper geodesic metric space $\XX$ is exponentially bounded if and only if it
is tempered.
\end{thm}

The proof consists in checking two following easy claims.

\begin{prop}
If a geodesic metric space $\XX$ carries a $(r_1,r_2)$-tempered measure with
$r_2/r_1>3$, then it is exponentially $r_1$-bounded.
\end{prop}

\begin{proof} \label{pr:tempered->bounded}
For simplicity put $\rho=r_1$ and $\e=r_2-3r_1$. For any $r\ge\rho$ let $N_r$
be the supremum (over $x\in\XX$) of cardinalities of disjoint systems of
$\rho$-balls in the balls $B(x,r)$. We shall show that the growth of $N_r$ is
at most exponential.

Let $\{x_i\}_{i=1}^N$ be a maximal system of points in the ball $B(x,r)$ such
that the balls $B(x_i,\rho)$ are pairwise disjoint. Since this system is
maximal, for any point $x'\in B(x,r-\rho)$ there is a point $x_i$ with
$d(x',x_i)\le 2\rho$, which means that $2\rho$-balls centered at the points
$x_i$ cover the ball $B(x,r-\rho)$, and, therefore, $r_2$-balls centered at
$x_i$ cover the ball $B(x,r+\e)$ for any $\e>0$. Thus,
$$
m B(x,r+\e) \le N \ov C(r_2) \le N_r \ov C(r_2) \;.
$$

On the other hand, the cardinality of any disjoint system of $\rho$-balls in
$B(x,r+\e)$ obviously does not exceed the ratio $m B(x,r+\e)/\un C(\rho)$,
whence
$$
\frac{N_{r+\e}}{N_r} \le \frac{\ov C(r_2)}{\un C(r_1)} \;,
$$
which implies the claim.
\end{proof}

\begin{prop}[{cf. \cite[Section 6.2]{Adams96}}] \label{pr:bounded->tempered}
If a geodesic metric space $\XX$ is exponentially $\rho$-bounded, then it
carries a $(9\rho,\infty)$-tempered measure. Moreover, if $\XX$ is proper, then
such a measure can be chosen to be $\Iso(\XX)$-invariant.
\end{prop}

\begin{proof}
Without the requirement of $\Iso(\XX)$-invariance the measure in question can
be easily taken to be, for instance, the sum of $\d$-measures on a maximal
$3\rho$-separated system of points in $\XX$. However, in order to obtain an
$\Iso(\XX)$-invariant measure a little bit more work has to be done.

Since the space $\XX$ is proper, the group $\Iso(\XX)$ is locally compact with
respect to the pointwise convergence topology. For any orbit
$\OO=\Iso(\XX)o,\,o\in\XX,$ choose the left-invariant Haar measure $\m_\OO$ on
$\Iso(\XX)$ normalized by the condition
\begin{equation} \label{eq:normalization}
\m_\OO\{g: d(o,go)\le 6\rho\}=1 \;,
\end{equation}
and denote by $m_\OO$ the image of $\m_\OO$ under the map $g\mapsto go$.
Obviously, the measure $m_\OO$ is $\Iso(\XX)$-invariant, it does not depend on
the choice of the point $o$ from the orbit $\OO$, and
$$
m_\OO B(x,6\rho)=1 \qquad\forall\,x\in\OO \;.
$$

Now take a maximal system of $\Iso(\XX$)-orbits $\OO_i$ in $\XX$ with the
property that the pairwise distances between different orbits from this system
are at least $3\rho$. The measure $m = \sum_i m_{\OO_i}$ is clearly
$\Iso(\XX)$-invariant. We claim that it is $(9\rho,\infty)$-tempered.

First of all, since the system $\{\OO_i\}$ is maximal, any point $x\in\XX$ lies
within the distance $3\rho$ from at least one of the orbits $\OO_i$. Therefore,
the ball $B(x,9\rho)$ contains a $6\rho$-ball centered at a point $x'\in\OO_i$,
so that
$$
m B(x,9\rho) \ge m_{\OO_i} B(x',6\rho) \ge 1
$$
by the definition of the measures $m_{\OO_i}$.

In order to obtain an upper estimate for the measures $m B(x,r)$ first notice
that the exponential $\rho$-boundedness of the space $\XX$ implies that the
ball $B(x,r)$ can be covered by at most $e^{ar}$ balls $B(x_i,3\rho)$ (see the
beginning of the proof of \propref{pr:tempered->bounded}). Thus, for any
$\Iso(\XX)$-orbit $\OO$ the intersection $B(x,r)\cap\OO$ can be covered by at
most $e^{ar}$ balls of radius $6\rho$ centered at points from $\OO$. Therefore,
\begin{equation} \label{eq:ear}
m_\OO B(x,r) \le e^{ar} \;.
\end{equation}

It remains to notice that the ball $B(x,r)$ can only intersect a uniformly
bounded number of orbits $\OO_i$. Indeed, choose a point $x_i$ on each orbit
$\OO_i$ intersecting $B(x,r)$. Then the pairwise distances between the points
$x_i$ are at least $3\rho$ by the definition of the system $\{\OO_i\}$.
Therefore, the $\rho$-balls centered at the points $x_i$ are all pairwise
disjoint and contained in $B(x,r+\rho)$, so that by the exponential
$\rho$-boundedness of the space $\XX$ the number of these balls (and therefore
of the orbits $\OO_i$ intersecting $B(x,r)$) does not exceed $e^{a(r+\rho)}$.
Combining this estimate with the inequality \eqref{eq:ear} gives
$$
m B(x,r) \le e^{a(2r+\rho)} \;.
$$
\end{proof}

Another bounded geometry type condition which can be imposed on a group of
isometries of a metric space is that of \emph{finiteness of the critical
exponent}.

\begin{defn}[{\cite{Burger-Mozes96}}]
Given a positive Radon measure $m$ on a metric space $\XX$, the number
(possibly infinite!)
\begin{equation} \label{eq:def-cr}
\d_{cr}(m) = \inf \left\{ \d\ge 0 : \int_\XX e^{- \d d(o,x)}\, dm(x) < \infty
\right\}
\end{equation}
is independent of the choice of the reference point $o\in\XX$, and is called
the \emph{critical exponent} of the measure $m$. The critical exponent
$\d_{cr}(G)$ of a closed subgroup $G\subset\Iso(\XX)$ is by definition the
critical exponent of a positive $G$-invariant measure supported on a $G$-orbit
in $\XX$. Equivalently, for any two points $o,x\in\XX$
$$
\d_{cr}(G) = \inf \left\{ \d\ge 0 : \int_G e^{- \d d(o,gx)}\, d\m(g) < \infty
\right\} \;,
$$
where $\m$ is a left-invariant Haar measure on $G$.
\end{defn}

Generally speaking, finiteness of the critical exponent of the group
$\Iso(\XX)$ does not imply exponential boundedness of the space $\XX$, because
the orbits of $\Iso(\XX)$ may be ``too small'' to provide any information about
the whole space $\XX$. However, the converse implication is easily seen to be
true:

\begin{prop} \label{pr:bounded->exponent}
If the space $\XX$ is exponentially bounded, then the critical exponent of any
closed subgroup $G\subset\Iso(\XX)$ is finite.
\end{prop}

\begin{proof}
Fix a $G$-invariant measure $m$ on the $G$-orbit of a reference point
$o\in\XX$. By changing variables the integral \eqref{eq:def-cr} from the
definition of the critical exponent can be presented as
\begin{equation} \label{eq:exponent}
\begin{aligned}
\int_\XX e^{- \d d(o,x)}\, dm(x)
 &= \int_0^1 m \left\{x: e^{- \d d(o,x)}\ge t\right\}\,dt \\
 &= \d \int_0^\infty m B(o,r) \, e^{-\d r}\,dr \;.
\end{aligned}
\end{equation}
Since exponential boundedness of $\XX$ implies that $m B(o,r)$ is dominated by
an exponential function of $r$ (cf. the proof of
\propref{pr:bounded->tempered}), the integral \eqref{eq:exponent} is finite if
$\d$ is big enough.
\end{proof}

\begin{rem}
In fact, it would be sufficient to look just at the full group $\Iso(\XX)$,
because for any closed subgroup its critical exponent does not exceed that of
the group $\Iso(\XX)$, see \cite[Lemma 1.5]{Burger-Mozes96}.
\end{rem}

\begin{rem}
If the group $\Iso(\XX)$ acts on $\XX$ cocompactly, then $\XX$ is easily seen
to be exponentially bounded (because the exponential rate of growth of the
locally compact group $\Iso(\XX)$ is finite), which by
\propref{pr:bounded->exponent} implies finiteness of the critical exponent of
$\Iso(\XX)$ (cf. \cite[Proposition 1.7]{Burger-Mozes96}).
\end{rem}

\subsection{Boundary maps}

\begin{thm} \label{th:lambda}
For any exponentially bounded hyperbolic space $\XX$ there exists a sequence of
$\Iso(\XX)$-equivariant Borel maps $\l_n:\XX\times\part \XX\to\Pr(\XX)$ such
that
\begin{equation} \label{eq:ass}
\|\l_n(x,\g)-\l_n(x',\g)\|\to 0 \qquad \forall\, x,x'\in \XX,\,\g\in\part \XX
\;.
\end{equation}
\end{thm}

For proving \thmref{th:lambda} we shall need the following simple inequality:

\begin{lem} \label{lem:sandwich}
Let $\{Z_k\}_{k=1}^\infty,\{Z'_k\}_{k=1}^\infty$ be two increasing sequences of
measurable sets of finite positive measure in a measure space $(\XX,m)$ such
that for a certain integer $\tau>0$
\begin{equation} \label{eq:sandwich}
 Z_k \subset Z'_{k+\tau} \;,
 \qquad Z'_k \subset Z_{k+\tau} \qquad\forall\, k\ge 1 \;,
\end{equation}
and let
$$
\l_n=\frac1n\sum_{k=1}^nm_{Z_k} \;, \qquad \l'_n=\frac1n\sum_{k=1}^nm_{Z'_k}
$$
be the Cesaro averages of the sequences $\{m_{Z_k}\},\{m_{Z'_k}\}$,
respectively. Then
$$
\| \l_n - \l'_n \| \le \frac{2\tau}n + \frac{4(n-\tau)}n
 \left[1- \left(\frac{m Z_1}{m Z_{n+\tau}}\right)^{\frac{2\tau}{n-\tau}}
 \right]  \qquad \forall\, n>\tau \;.
$$
\end{lem}

\begin{proof}
First notice that for any $k>\tau$
$$
\| m_{Z_{k-\tau}} - m_{Z'_k} \|
 = 2\left(1-\frac{m Z_{k-\tau}}{m Z'_k}\right)
 \le 2\left(1-\frac{m Z_{k-\tau}}{m Z_{k+\tau}}\right) \;,
$$
and the same inequality holds for $\| m_{Z_{k-\tau}} - m_{Z_k} \|$, so that
$$
\| m_{Z_k} - m_{Z'_k} \|
 \le 4\left(1-\frac{m Z_{k-\tau}}{m Z_{k+\tau}}\right) \qquad\forall\,k>\tau\;.
$$
Then
$$
\begin{aligned}
 \| \l_n - \l'_n \|
 \le \frac1n \sum_{k=1}^n \| m_{Z_k} -m_{Z'_k} \|
 &\le \frac{2\tau}n
 + \frac4n \sum_{k=\tau+1}^n \left(1- \frac{m Z_{k-\tau}}{m Z_{k+\tau}} \right) \\
 &\le \frac{2\tau}n
 + \frac{4(n-\tau)}n\left[1 - \left(\prod_{k=\tau+1}^n
 \frac{m Z_{k-\tau}}{m Z_{k+\tau}}\right)^{\frac1{n-\tau}}\right] \;,
\end{aligned}
$$
whence the claim.
\end{proof}

\begin{cor}
If in \lemref{lem:sandwich} the growth of the sets $Z_n$ with respect to the
measure $m$ is subexponential, i.e., $\bigl(m Z_n\bigr)^{1/n}\to 1$, then $\|
\l_n - \l'_n\|\to 0$.
\end{cor}

\begin{proof}[Proof of \thmref{th:lambda}]
Given points $x\in \XX,\g\in\part \XX$ and positive integers $n\ge k\ge 1$ put
$$
Y(x,\g,n,k) = \bigl\{\xi(n): \xi\in\RR, \, d(x,\xi(0))\le k,\, \xi(\infty)=\g
\bigr\} \;,
$$
where $\RR$ is the set of all geodesic rays in $\XX$. We shall fix a
$(r,\infty)$-tempered $\Iso(\XX)$-invariant measure $m$ on $\XX$ provided by
\thmref{thm:measure}, and denote by $Z(x,\g,n,k)$ the closed $r$-neighbourhood
of $Y(x,\g,n,k)$.

As it follows from \propref{prop:sausage}, for any ray $\eta$ joining $x$ and
$\g$ the set $Z(x,\g,n,k)$ is contained in the closed $(C_1+r)$-neighbourhood
of the geodesic segment $\eta([n-k,n+k])$. Therefore, since the measure $m$ is
$(r,\infty)$-tempered, the values $m Z(x,\g,n,k)$ are bounded away from $0$ and
uniformly bounded from above by an affine function of $k$.

For any $\tau\ge d(x,x')$ the finite sequences of sets
$$
\bigl\{Z(x,\g,n,k)\bigr\}_{k=1}^n \qquad\text{and}\qquad
\bigl\{Z(x',\g,n,k)\bigr\}_{k=1}^n
$$
are $\tau$-sandwiched
in the sense of formula \eqref{eq:sandwich}, whence, by putting
$$
\l_n(x,\g) = \frac1n \sum_{k=1}^n m_{Z(x,\g,n,k)}
$$
and applying \lemref{lem:sandwich} we obtain the claim.
\end{proof}

\begin{thm} \label{th:lambda-CAT}
If $G$ is a closed group of isometries of a \emph{CAT$(-1)$} space $\XX$ with
finite critical exponent, then there exists a sequence of $G$-equivariant Borel
maps $\l_n:\XX\times\part \XX\to\Pr(\XX)$ such that
\begin{equation} \label{eq:ass-CAT}
\|\l_n(x,\g)-\l_n(x',\g)\|\to 0 \qquad \forall\, x,x'\in \XX,\,\g\in\part \XX
\;.
\end{equation}
\end{thm}

\begin{proof}
Let us fix a number $\d>\d_{cr}$, a reference point $o\in\XX$ and a
$G$-invariant Radon measure $m$ on the orbit $Go$. For any $x\in\XX$ denote by
$\nu_x=\nu_x^\d$ the probability measure on the orbit $Go$ with the density
$$
\frac{e^{-\d d(x,y)}}{\int_\XX e^{-\d d(x,y)} dm(y)}
$$
with respect to the measure $m$. We shall call the measures $\nu_x$
\emph{pre-Patterson} as their limits when $\d$ tends to the critical exponent
$\d_{cr}$ provide the famous Patterson measures on the boundary $\part\XX$, see
\cite{Patterson76}, \cite{Burger-Mozes96}. Obviously, the map $x\mapsto\nu_x$
is $G$-equivariant. This fact, together with the local Lipschitz continuity of
the map $x\mapsto\nu_x$ and boundary convergence of geodesics (guaranteed by
the CAT($-1$) condition) yields the proof.

Since for any $x,x',y\in\XX$
$$
{e^{-\d d(x,x')}}
 \le \frac{e^{-\d d(x',y)}}{e^{-\d d(x,y)}}
 \le {e^{\d d(x,x')}} \;,
$$
the definition of the pre-Patterson measures implies that
\begin{equation} \label{eq:close}
\| \nu_x - \nu_{x'} \| \le 3\d\, d(x,x')
\end{equation}
provided $d(x,x')$ is small enough.

Now put
$$
\l_n(x,\g) = \frac1n \int_0^n \nu_{\xi(t)} \,dt \;,
$$
where $\xi$ is the geodesic ray joining $x$ and $\g$. By the CAT($-1$) property
for any other point $x'\in\XX$
$$
d(\xi(t),\xi'(t+T)) \toto_{t\to\infty} 0 \;,
$$
where $\xi'$ is the geodesic ray joining $x'$ and $\g$, and $T$ is the value of
the Busemann cocycle $\b_\g(x,x')$ (see \remref{rem:Busemann}), which in
combination with formula \eqref{eq:close} implies the claim.
\end{proof}

\begin{rem} \label{rem:uniform}
It can be easily seen from the proofs of \thmref{th:lambda} and
\thmref{th:lambda-CAT} that for any compact subset $K\subset\XX$ the
convergence in formulas \eqref{eq:ass} and \eqref{eq:ass-CAT}, respectively, is
uniform on $K\times K\times\part\XX$.
\end{rem}

\section{Hyperbolic equivalence relations} \label{sec:relations}

\subsection{Graphed equivalence relations} \label{sec:equiv relations}

We begin with recalling the basic notions of the theory of discrete equivalence
relations, see \cite{Feldman-Moore77}. Let $X$ be a standard Borel space. An
equivalence relation $R$ on $X$ is called \emph{standard} if it is a Borel
subset of $X\times X$, and it is called \emph{countable} if the equivalence
class (the \emph{leaf}) $[x]=R(x)=\{y: (x,y)\in R\}$ of any point $x\in X$ is
at most countable. A standard countable equivalence relation is also called
\emph{discrete}.

\begin{prop}[{cf. \cite[proof of Theorem 1]{Feldman-Moore77}}] \label{pr:alpha}
For any discrete equivalence relation $R$ with infinite equivalence classes
there exists a sequence $\{\a_n\}_{n=1}^\infty$ of measurable maps $\a_n:X\to
X$ such that their graphs are pairwise disjoint and their union is $R$. In
other words, for any $x\in X$ the sequence $\a_n(x)$ determines an ordering of
the class $[x]$, and these orderings depend on $x$ measurably.
\end{prop}

A standard equivalence relation $R$ is called \emph{non-singular} with respect
to a Borel probability measure $\mu$ on $X$ (or, equivalently, the measure
$\mu$ is \emph{quasi-invariant} with respect to $R$) if for any subset
$A\subset X$ with $\mu(A)=0$ its \emph{saturation} $[A]=\bigcup_{x\in A}[x]$
also has measure 0. A non-singular equivalence relation $(X,\mu,R)$ is
\emph{ergodic} if there are no non-trivial saturated sets, it is
\emph{conservative} if all ergodic components are uncountable, and it is
\emph{dissipative} (or of type I) if any ergodic component consists of a single
equivalence class, or, equivalently, there exists a leafwise constant
measurable map $\f:X\to X$ with $(x,\f(x))\in R$ for a.e.~$x\in X$. In
particular, dissipative ergodic equivalence relations are precisely those which
consist of a single class. Any non-singular equivalence relation can be
uniquely decomposed into a disjoint union of its conservative and dissipative
parts.

Integrating the counting measures on the fibres of the left $(x,y)\mapsto x$
and the right $(x,y)\mapsto y$ projections from $R$ onto $X$ by the measure
$\mu$ gives the {\it left\/} $d\Mu(x,y)=d\mu(x)$ and the {\it right\/}
$d\check\Mu(x,y)=d\Mu(y,x)=d\mu(y)$ {\it counting measures\/} on $R$,
respectively. The measures $\Mu$ and $\check\Mu$ are equivalent iff $\mu$ is
quasi-invariant. If $\Mu=\check\Mu$, then the measure $\mu$ is called
$R$-invariant.

A (non-oriented) \emph{graph structure} on an equivalence relation $(X,\mu,R)$
is determined by a measurable subset $K\subset R$ which is symmetric and does
not intersect the diagonal $\{(x,x)\}\subset R$. Two points $x,y\in X$ are then
joined with an edge iff $(x,y)\in K$. We shall call $(X,\mu,R,K)$ a
\emph{graphed equivalence relation} \cite{Adams90}. Denote by $[x]_K$ the
equivalence class $[x]$ endowed with the graph structure $K$.

\begin{defn}
A graphed equivalence relation $(X,\mu,R,K)$ is \emph{hyperbolic} if a.e.~
leafwise graph $[x]_K$ is a hyperbolic graph (i.e., is Gromov hyperbolic with
respect to the graph metric assigning length 1 to all edges). By $\ov{[x]}_K$
and $\part [x]_K$ we denote the hyperbolic compactification of the graph
$[x]_K$ and its boundary, respectively. As usual, $C\bigl(\ov{[x]}_K\bigr)$
denotes the space of continuous real valued functions on $\ov{[x]}_K$. Below we
shall always deal with a fixed graph structure $K$ and usually omit the
corresponding subscript.
\end{defn}

\begin{rem}
The function assigning to any point $x\in X$ the hyperbolicity constant of the
graph $[x]$ is easily seen to be measurable (because the hyperbolicity constant
is defined in terms of geodesic triangles in $[x]$). Being leafwise constant,
this function must therefore be constant on the ergodic components of $R$.
\end{rem}

Examples of hyperbolic equivalence relations include, in particular,
\begin{enumerate}
\item
the leafwise Cayley graphs of a free action of a word hyperbolic group, see
\exref{ex:orbits} for more details;
\item
equivalence relations on transversals to hyperbolic foliations, see below
\secref{sec:fol};
\item
\emph{treed equivalence relations} (for which a.e.~leafwise graph is a tree),
see \cite{Adams90};
\item
equivalence relations associated with certain fractal sets, see
\cite{Kaimanovich02}.
\end{enumerate}

In all these examples the hyperbolic boundaries of different equivalence
classes can be naturally identified with (a subspace of) the same topological
space (see the discussion of these examples at the end of
\secref{sec:measurable systems}). However, there is no reason for this to be
the case in general (it would be interesting to have an explicit example of
this kind), which is why we shall need a technique allowing one to deal with
the general situation.

\subsection{The bundle of leafwise continuous functions} \label{sec:continuous
bundle}

Recall that a \emph{measurable structure} on a bundle of Banach spaces
$\{E_x\}$ over a Lebesgue ($\equiv$ standard measure) space $(X,\mu)$ is a
family $\M$ of sections (called \emph{measurable}) $\s:x\mapsto \s_x\in E_x$
such that
\begin{itemize}
  \item[(i)]
the function $x\mapsto \|\s_x\|$ is measurable for any section $\s\in\M$;
  \item[(ii)]
the family $\M$ is closed under addition, multiplication by measurable scalar
functions and passing to pointwise limits in the norm topology.
\end{itemize}
The couple $(\{E_x\},\M)$ is then called a \emph{measurable bundle of Banach
spaces} over $(X,\mu)$, see \cite[Chapter II.4]{Fell-Doran88}, \cite[Appendix
A.3]{Anantharaman-Renault00}. A measurable bundle \linebreak $(\{E_x\},\M)$ is
\emph{separable} if there exists a sequence of sections $\s^n\in\M$ such that
for a.e.~$x\in X$ the sequence $\{\s^n_x\}$ is dense in $E_x$.

\begin{prop} \label{pr:separab}
Given a hyperbolic equivalence relation $(X,\mu,R,K)$, denote by
$\M=\M(X,\mu,R,K)$ the space of all measurable functions $F:R\to\R$ such that
for a.e.~$x\in X$ the function $F(x,\cdot)$ extends to a continuous function on
the completion $\ov{[x]}$. Then $\M$ is a measurable structure on the bundle of
Banach spaces $\bigl\{C\bigl(\ov{[x]}\bigr)\bigr\},\,x\in X$, and the
measurable bundle $\bigl(\bigr\{C\bigl(\ov{[x]}\bigl)\bigr\},\M\bigr)$ is
separable.
\end{prop}

\begin{proof}
Verification of the conditions from the definition of a measurable structure
being straightforward, we only have to check the separability. Since the
leafwise hyperbolic compactifications $\ov{[x]}$ are completions of the
equivalence classes $[x]$ with respect to the metrics $\rho_x$ \eqref{eq:rho},
and these metrics depend on $x$ measurably, the functions
$$
F_n(x,y)=\rho_x(y,\a_n(x)) \;,
$$
where $\a_n$ are the maps from \propref{pr:alpha}, belong to $\M$. For any
$x\in X$ the family of functions $\{F_n(x,\cdot)\}$ separates points of
$\ov{[x]}$. By adding to this family the constant function and applying the
Weierstra{\ss}--Stone theorem we obtain the claim.
\end{proof}

The \emph{dual bundle} $(\{E^*_x\},\M^*)$ of a separable measurable bundle
$(\{E_x\},\M)$ is the bundle of dual Banach spaces $E^*_x$ endowed with the
\emph{weak$^*$ measurable structure} $\M^*$ which consists of all sections
$\s^*:x\mapsto \s^*_x\in E^*_x$ such that the map $x\mapsto \langle \s_x,\s_x^*
\rangle$ is measurable for any section $\s\in\M$, see
\cite[Lemma~A.3.7]{Anantharaman-Renault00}.

\medskip

Returning to the measurable bundle
$\bigl(\bigl\{C\bigl(\ov{[x]}\bigr)\bigr\},\M\bigr)$ of leafwise continuous
functions, which is associated with the hyperbolic equivalence relation
$(X,\mu,R,K)$, notice that for any space $C\bigl(\ov{[x]}\bigr)$ its dual is
the space $M\bigl(\ov{[x]}\bigr)$ of signed Borel measures (equipped with the
total variation norm) on the completion $\ov{[x]}$. Therefore, in view of
\propref{pr:separab} we have a measurable structure $\M^*$ on the bundle
$\bigl\{M\bigl(\ov{[x]}\bigr)\bigr\}$. Denote by $\M^*_\circ$ and $\M^*_\part$
the measurable structures induced by $\M^*$ on the \emph{interior subbundle}
$\{M([x])\}$ and the \emph{boundary subbundle} $\{M(\part[x])\}$, respectively.

\begin{rem} \label{rem:Hahn}
One can easily check that the splitting of the bundle
$\bigl\{M\bigl(\ov{[x]}\bigr)\bigr\}$ into the direct sum of the interior and
boundary subbundles is measurable in the sense that for any measurable section
of $\bigl\{M\bigl(\ov{[x]}\bigr)\bigr\}$ its interior and boundary parts are
measurable too. Also, the \emph{Hahn decomposition} of any measurable section
of $\bigl\{M\bigl(\ov{[x]}\bigr)\bigr\}$ into the difference of its positive
and negative parts is measurable.
\end{rem}

\begin{defn} \label{def:measure systems}
Measurable sections of the bundle $(\{M([x])\},\M^*_\circ)$ (resp., of the
bundle $(\{M(\part [x])\},\M^*_\part)$) are called \emph{measurable systems of
interior $($resp., \linebreak boundary$)$ measures} on a hyperbolic equivalence
relation $(X,\mu,R,K)$.
\end{defn}

\subsection{Measurable systems of boundary measures} \label{sec:measurable
systems}

Since the functions
$$
\D_n(x,y) =
  \begin{cases}
    1 \;, & y=\a_n(x) \;, \\
    0 \;,& \text{otherwise} \;,
  \end{cases}
$$
where $\a_n$ are the maps from \propref{pr:alpha}, obviously belong to $\M$,
the measurable structure $\M^*_\circ$ of the interior bundle $\{M([x])\}$
consists precisely of those systems of measures $\{\l_x\}$ for which $\l_x(y)$
is measurable as a function on $R$. We shall now give an explicit description
of the measurable structure $\M^*_\part$ of the boundary bundle
$\{M(\part[x])\}$.

\begin{thm} \label{th:continuous}
The measurable structure $\M^*_\part$ of the boundary bundle consists of
pointwise weak$^*$ limits $($in the spaces $M\bigl(\ov{[x]}\bigr),\,x\in X)$ of
sequences of measurable sections of the interior bundle $\{M([x])\}$.
\end{thm}

\begin{proof}
The definition of the dual bundle implies that its measurable structure $\M^*$
is closed with respect to passing to pointwise weak$^*$ limits, so that we only
have to prove that any measurable system of boundary measures $\l\in\M^*_\part$
can be approximated by a sequence of measurable systems of interior measures
$\l^n\in\M^*_\circ$. In other words, for any boundary measure $\l_x\in
M(\part[x])$ we have to construct a sequence of interior measures $\l^n_x\in
M([x])$ converging to $\l_x$, and to do it in a measurable way with respect to
$x$.

Let us begin with the following observation. Any finite continuous partition of
unity
$$
\1=\sum_{i\in I}\f_i
$$
on a compact space $K$ determines a norm 1 linear map $\pi$ from the space
$M(K)$ of (signed) Borel measures on $K$ to the space $M(I)$ of (signed)
measures on the index set $I$ by the formula
\begin{equation} \label{eq:partition}
\pi\th(i) = \langle \f_i, \th \rangle \;.
\end{equation}
Suppose now that for any $x\in X$ and any integer $n$ we have a partition of
unity on the completion $\ov{[x]}$ parameterized by the points from the
$n$-sphere of the leafwise graph metric $S(x,n)\subset [x]$. These partitions
then determine a sequence of maps $\pi_n$ from $M\bigl(\ov{[x]}\bigr)$ to
$M(S(x,n))$. We shall show that these partitions of unity can be constructed
measurably with respect to $x$, and in such a way that $\pi_n\th\to\th$ for any
boundary measure $\th\in M(\part[x])$.

\medskip

For the moment we fix a point $x\in X$, an integer $n>0$, and a real $\e>0$.
Recall that the hyperbolic compactification of the equivalence class $[x]$ is
its completion with respect to the metric $\rho_x$ \eqref{eq:rho}. Now for any
$s\in S=S(x,n)$ define the non-negative function $f_s\in C\bigl(\ov{[x]}\bigr)$
by putting for $z\in\ov{[x]}$
$$
f_s(z) =
  \begin{cases}
    1 - \displaystyle{\frac{\rho_x(z,s) - \rho_x(z,S)}\e} &, \quad \rho_x(z,s) \le
    \rho_x(z,S) + \e \;; \\
    0 &, \quad \text{otherwise} \;.
  \end{cases}
$$
Then the functions
\begin{equation} \label{eq:defphi}
\f_s = \frac{f_s}{\sum_{s\in S} f_s}
\end{equation}
form a partition of unity on $\ov{[x]}$ parameterized by the points $s\in S$,
and these functions have the property that
\begin{equation} \label{eq:vanish}
\f_s(z)=0 \qquad \forall\, z\in\ov{[x]}: \rho_x(z,s)\ge\rho_x(z,S)+\e \;.
\end{equation}

As it follows from the definition of the metric $\rho_x$ \eqref{eq:rho},
\begin{equation} \label{eq:ean}
\rho_x(\g,S) = e^{-an} \qquad \forall\,\g\in\part[x] \;,
\end{equation}
where $a=1/15\d$. Indeed, $\rho_x(y,s)\ge e^{-an}$ for any two points $y\in
[x]$ and $s\in S$, so that $\rho_x(\g,s)\ge e^{-an}$. On the other hand, if
$\xi$ is a geodesic ray joining $x$ and $\g$, then
$$
 \rho_x\bigl(\g,\xi(n)\bigr)
 = \lim_{m\to\infty}\rho_x\bigl(\xi(m),\xi(n)\bigr)
 = e^{-an} \;.
$$

From now on we shall assume, for the sake of concreteness, that in the above
construction
$$
\e=\e(n)=e^{-an} \;.
$$
Then formulas \eqref{eq:vanish} and \eqref{eq:ean} imply that
\begin{equation} \label{eq:2ean}
\f_s(\g) = 0 \qquad \forall\, \g\in\part[x]: \rho_x(\g,s)\ge 2e^{-an} \;.
\end{equation}

Denote now by
$$
\pi_n:M\bigl(\ov{[x]}\bigr)\to M(S)
$$
the map \eqref{eq:partition} from the space of measures on $\ov{[x]}$ to the
space of measures on $S=S(x,n)$ determined by the partition of unity
$\{\f_s\}_{s\in S}$ \eqref{eq:defphi}. Then by \eqref{eq:2ean}
\begin{equation} \label{eq:supp}
\supp\pi_n\d_\g \subset S \cap B_x(\g,2e^{-an}) \qquad \forall\,\g\in\part[x]
\;,
\end{equation}
and the linearity of the maps $\pi_n$ implies the weak$^*$ convergence of the
sequence $\pi_n\th$ to $\th$ for any boundary measure $\th\in M(\part[x])$.
Obviously, the whole construction is measurable in $x$.
\end{proof}

\begin{prop} \label{pr:atomic decomposition}
For any measurable system of boundary measures on a hyperbolic equivalence
relation its continuous and atomic parts are also measurable.
\end{prop}

\begin{proof}
In view of \remref{rem:Hahn} it is enough to establish the claim for systems of
non-negative measures only. For any such system $\l=\{\l_x\}$ it is sufficient
to show that the system $\wh\l=\{\wh\l_x\}$ which consists of maximal weight
atoms of the measures $\l_x$ is also measurable.

We shall use the approximation of the measures $\l_x$ by the interior measures
$\pi_n\l_x$ constructed in the proof of \thmref{th:continuous}. Let us first
notice that, as it follows from formulas \eqref{eq:ean} and \eqref{eq:supp},
for any point $\g\in\part[x]$ there exists $s\in S(x,n)$ such that
\begin{equation} \label{eq:supp3}
\supp \pi_n\d_\g \subset S(x,n) \cap B_x(s,3e^{-an}) \;.
\end{equation}
Formulas \eqref{eq:supp} and \eqref{eq:supp3} show that the maximal weight of
the measure $\l_x$ can be expressed in terms of the measures $\pi_n\l_x$ as
$$
W(x) = \lim_{n\to\infty} \max_{s\in S(x,n)} (\pi_n\l_x) \bigl[ B_x(s,3e^{-an})
\bigr] \;,
$$
and, in particular, the function $x\mapsto W(x)$ is measurable. Let us now
denote by $\wh\l^n_x$ the restriction of the measure $\pi_n\l_x$ to the
$3e^{-an}$-neighbourhood (in the metric $\rho_x$) of the set
$$
\Bigl\{s\in S(x,n): (\pi_n\l_x) \bigl[ B_x(s,3e^{-an}) \bigr] \ge W(x) \Bigr\}
\;.
$$
Then for any $x\in X$ the sequence $\wh\l^n_x$ converges (in the weak$^*$
topology) to the measure $\wh\l_x$. Since for any $n$ the system of measures
$\{\wh\l^n_x\}$ is measurable, by \thmref{th:continuous} the system $\wh\l_x$
is also measurable.
\end{proof}

In the situations when the hyperbolic boundaries of different equivalence
classes can be identified with (a subspace of) the same topological space, the
notion of measurability introduced above coincides with the ``usual''
measurability as one can easily verify by using \thmref{th:continuous}.

\begin{ex} \label{ex:orbits}
Recall that the classes of the \emph{orbit equivalence relation} $R_G$ of a
measure class preserving action of a countable group $G$ on a measure space
$(X,\mu)$ are the orbits of the action. Any symmetric generating set $A\subset
G$ determines the graph structure $K=\{(x,gx):x\in X,g\in A\}$ on $R_G$. If the
action of $G$ is free (mod 0), then the map $g\mapsto g^{-1}x$ establishes an
isomorphism of the (right) Cayley graph of $G$ with a.e.~graph $[x]$. If the
group $G$ is in addition \emph{word hyperbolic} (i.e., the Cayley graph $(G,A)$
is hyperbolic), then the orbit equivalence relation $R_G$ is obviously
hyperbolic, and the above isomorphism extends to a homeomorphism $\f_x$ from
the hyperbolic boundary $\part G$ of the group $G$ to the hyperbolic boundary
$\part[x]$ of the equivalence class $[x]$ for a.e.~point $x\in X$. Therefore,
in this situation $M(\part[x])\cong M(\part G)$ for a.e.~$x\in X$, so that the
boundary bundle is isomorphic to a constant Banach bundle, and
\thmref{th:continuous} implies that the measurability of a system of boundary
measures $\l_x\in M(\part[x])$ in the sense of
\defref{def:measure systems} is equivalent to the measurability of the map
$x\mapsto\f_x^{-1}(\l_x)$ from $X$ to $M(\part G)$.
\end{ex}

\begin{ex}
Any simply connected Riemannian manifold $M$ with pinched negative sectional
curvatures is $\d$-hyperbolic, and the exponential map at any point $x\in M$
establishes a homeomorphism between the unit tangent sphere at the point $x$
and the hyperbolic boundary of $M$ (in this situation the hyperbolic
compactification coincides with the visibility compactification). Therefore,
given a hyperbolic Riemannian measurable foliation, the leafwise hyperbolic
boundaries are homeomorphic to spheres, and the boundary bundle for the
equivalence relation on any uniform system of transversals (see
\secref{sec:fol}) is isomorphic to a constant Banach bundle. The measurability
of a system of boundary measures on this equivalence relation is then
equivalent to the measurability of the corresponding system of measures on the
leafwise unit tangent spheres of the foliation.
\end{ex}

\begin{ex}
Yet another example when the boundary bundle is isomorphic to a constant Banach
bundle (more precisely, to a subbundle of a constant bundle) is provided by
\emph{treed equivalence relations} (see \cite{Adams90}) for which a.e.~graph
$[x]$ is a tree. In this case any such tree can be embedded into the universal
countable tree. Once again, \thmref{th:continuous} implies that the
measurability of a system of boundary measures in the sense of
\defref{def:measure systems} is equivalent to the measurability of the
corresponding map from $X$ to the space of boundary measures of the universal
tree.
\end{ex}

\subsection{Invariant systems of boundary measures}

\begin{defn} \label{def:inv-system}
A measurable system $\l=\{\l_x\}$ of boundary measures on a hyperbolic
equivalence relation $(X,\mu,R,K)$ is called \emph{invariant} if it is a.e.~
constant, i.e., $\l_x=\l_y$ for $\Mu$-a.e.~$(x,y)\in R$.
\end{defn}

\begin{thm} \label{th:at most 2}
If $\l=\{\l_x\}$ is an invariant system of boundary measures on a conservative
hyperbolic equivalence relation $(X,\mu,R,K)$, then a.e.~measure $\l_x$ is
supported by at most 2 points.
\end{thm}

\begin{proof}
As it follows from the proof of \propref{pr:atomic decomposition}, the set $X'$
of all points $x\in X$ with $\card\supp\l_x\ge 3$ is measurable, so that
without loss of generality we may assume that $X'=X$. Further, again by the
proof of \propref{pr:atomic decomposition}, we may change the system $\l$ so
that the weights of atoms of a.e.~measure $\l_x$ be strictly less than
$\displaystyle\frac12$.

The claim is now a straightforward consequence of \propref{pr:bary}. The only
thing we have to take care of is the measurability of the quasi-barycenter set.
In order to check it we shall invoke \thmref{th:continuous} and approximate (in
the leafwise weak$^*$ topology) the system $\l$ by a sequence of measurable
systems of interior measures $\l^n=\{\l_x^n\}$. Define the integer valued
function $\f$ on $R$ as
$$
\f(x,y) = \limsup_{n\to\infty} \sum_z \b_z(x,y) \l_x^n(z) \;,
$$
where $\b_z$ is the distance cocycle \eqref{eq:bz}. Let $Y(x)\subset [x]$ be
the set of points where the function $\f(x,\cdot)$ attains its minimal value,
and put
$$
Z(x) = \bigcup_{x'\in[x]} Y(x') \;.
$$
Clearly, the set
$$
\bigl\{(x,y):x\in X,y\in Z(x)\bigr\} \subset R
$$
is measurable, and by \propref{pr:b-infinity} and \propref{pr:bary} $Z(x)$ is
contained in the finite set $\CC(\l,x,8C_2)$ for a.e.~$x\in X$. Thus, we have a
measurable map assigning to any point $x\in X$ a finite subset $Z(x)$ of the
equivalence class $[x]$, which contradicts the conservativity of $R$.
\end{proof}

\begin{defn} \label{def:elem}
A hyperbolic equivalence relation is \emph{elementary} if there exists a
measurable map assigning to a.e.~class $[x]$ a finite subset $A[x]\subset\part
[x]$, i.e., if the system of uniform probability measures on the sets $A[x]$ is
an invariant measurable system in the sense of \defref{def:measure systems}.
\end{defn}

In the setup of this definition the map $x\mapsto\card A[x]$ is measurable (see
the proof of \propref{pr:atomic decomposition}). Therefore, if an elementary
hyperbolic equivalence relation is ergodic, then almost all sets $A[x]$ must
have the same cardinality. If $q\ge 1$ is its value, then we shall call the
associated system of boundary measures and the equivalence relation
\emph{$q$-elementary}. Note that even if a hyperbolic equivalence relation is
elementary, a.e.~leafwise boundary may well be infinite (see the examples
below). \thmref{th:at most 2} implies

\begin{thm} \label{th:3types}
Let $(X,\mu,R,K)$ be an ergodic hyperbolic equivalence relation which has more
than one equivalence class. Then it admits an invariant measurable system of
boundary probability measures if and only if it belongs to one of the following
3 classes:
\begin{itemize}
\item[(i)]
$(X,\mu,R,K)$ admits a unique invariant measurable system of boundary measures
which is 1-elementary;
\item[(ii)]
$(X,\mu,R,K)$ admits a unique invariant measurable system of boundary measures
which is 2-elementary;
\item[(iii)]
$(X,\mu,R,K)$ admits two disjoint 1-elementary systems of boundary measures,
and any invariant system of boundary measures is their convex combination.
\end{itemize}
\end{thm}

\subsection{Examples of elementary hyperbolic equivalence relations}
\label{sec:examples}

We shall now give examples of hyperbolic equivalence relations from each of the
types described in \thmref{th:3types}. Take a \emph{non-elementary} word
hyperbolic group $G$, i.e., such that its hyperbolic boundary $\part G$
consists of at least 3 points (which implies that $\part G$ is uncountable, see
\cite{Gromov87}). For instance, one can take for $G$ a non-abelian free group.
Put
\begin{eqnarray*}
 & \part^2 G = \part G\times \part G \setminus\diag =
 \{(\g_1,\g_2):\g_i\in\part G,\g_1\neq\g_2\} \;, & \\
 & \part^3 G = \{(\g_1,\g_2,\g_3):\g_i\in\part G,\g_i\neq\g_j\} \;, &
\end{eqnarray*}
and, finally, let $\wt{\part^2 G}$ be the \emph{symmetrization} of $\part^2 G$,
i.e., the set of unordered pairs of elements $\g_1\neq\g_2\in\part G$.

\begin{prop} \label{pr:elem}
The orbit equivalence relation $R_G$ of a free measure class preserving action
of a word hyperbolic group $G$ on a measure space $(X,\mu)$ is 1-elementary
$($resp., 2-elementary$)$ if and only if there exists a measurable
$G$-equi\-va\-ri\-ant map $\pi:X\to\part G$ $($resp., $\pi:X\to \wt{\part^2
G})$.
\end{prop}

\begin{proof}
First notice that the homeomorphisms $\f_x:\part G\mapsto\part[x]$ (see
\exref{ex:orbits}) have the property that $\f_{gx}(g\g) = \f_x(\g)$ for any
$g\in G,\g\in\part G$ and a.e.~$x\in X$. Then equivariance of $\pi$ implies
that the map $x\mapsto\f_x\circ\pi(x)\in\part [x]$ is leafwise constant, i.e.,
$R_G$ is 1-elementary. Conversely, if $R_G$ is 1-elementary, and
$x\mapsto\F(x)\in\part[x]$ is the corresponding leafwise constant map, then the
map $\pi(x)=\f_x^{-1}\circ\F(x)$ is equivariant. With obvious modifications the
same argument works in the 2-elementary case as well.
\end{proof}

Let us fix a quasi-invariant probability measure $\mu$ on $\part G$ such that
the diagonal action of $G$ on $(\part^2 G,\mu\otimes\mu)$ is ergodic (for
instance, the harmonic measure of any non-degenerate symmetric random walk on
$G$ has this property, see \cite{Kaimanovich94}). On the other hand,
\propref{pr:bary} implies that the action of $G$ on $\part^3 G$ is proper (cf.
\cite[8.2.K]{Gromov87}). Therefore, this action is dissipative with respect to
any quasi-invariant measure.

\begin{ex}
\emph{A hyperbolic equivalence relation with a unique 1-elementary system of
boundary measures.} This is the orbit equivalence relation of the action of $G$
on $(\part G,\mu)$.

\propref{pr:elem} applied to the identity map on $\part G$ at once provides us
with a 1-elementary system. In order to prove its uniqueness we have, in view
of \thmref{th:3types}, to exclude existence of another 1-elementary system of
boundary measures disjoint with the one we have just constructed. By
\propref{pr:elem} existence of such a system is equivalent to existence of an
equivariant map $\pi:\part G\to\part G$ such that a.e.~$\pi(\g)\neq \g$. Such a
map $\pi$ would allow one to equivariantly embed $\part^2 G$ into $\part^3 G$
by the formula $(\g_1,\g_2)\mapsto (\g_1,\pi(\g_1),\g_2)$, which is impossible
because the action of $G$ on the space $(\part^2 G,\mu\otimes\mu)$ is ergodic,
whereas its action on $\part^3 G$ (with respect to the image of the measure
$\mu\otimes\mu$) must be dissipative.
\end{ex}

\begin{ex}
\emph{A hyperbolic equivalence relation with a unique system of \linebreak
boundary measures which is 2-elementary.} This is the orbit equivalence
relation of the action of $G$ on $\wt{\part^2 G}$ (endowed with the image of
the measure $\mu\otimes\mu$).

The system in question by \propref{pr:elem} corresponds to the identity map on
$\wt{\part^2 G}$. By \thmref{th:3types} in order to prove its uniqueness we
have to exclude the possibility that this system is the half-sum of two
1-elementary systems, i.e., in view of \propref{pr:elem}, that there exists an
equivariant map $\pi:\wt{\part^2 G}\to\part G$ with
$\pi(\{\g_1,\g_2\})\subset\{\g_1,\g_2\}$. Such a map $\pi$ would give a
$G$-equivariant measurable section of the projection $\part^2 G\to\wt{\part^2
G}$, which is impossible by the ergodicity of the action of $G$ on $(\part^2
G,\mu\otimes\mu)$.
\end{ex}

\begin{ex}
\emph{A hyperbolic equivalence relation with two 1-elementary systems of
boundary measures.} This is the orbit equivalence relation of the action of $G$
on $(\part^2 G,\mu\otimes\mu)$.
\end{ex}

\subsection{Foliations and equivalence relations} \label{sec:fol}

Recall that a \emph{foliation} $\FF$ is a smooth manifold decomposed into
smooth immersed submanifolds \emph{(leaves)} organized in a local product
structure (we shall use the same notation $\FF$ both for the foliation and for
its underlying manifold). In other words, $\FF$ can be covered by \emph{flow
boxes} $U\cong B\times T$ (where $B$ and $T$ are open Euclidean balls) endowed
with an additional product structure preserved by the transition maps. The sets
$\{x\}\times T,\,x\in B$ and $B\times\{t\},\,t\in T$ are called
\emph{transversals} and \emph{plaques}, respectively. By $\FF(x)$ we denote the
\emph{leaf} of the foliation $\FF$ passing through a point $x$.
\emph{Laminations} are defined in precisely the same way except for allowing
the transverse structure to be just continuous (rather than smooth), so that
for laminations the transversals $T$ are just topological spaces (without any
smooth structure). However, the difference between foliations and laminations
is insignificant for us, as we shall always (unless otherwise specified) be
working in the \emph{measure category}, and only require measurability of the
transverse structure of $\FF$ with respect to a certain \emph{quasi-invariant
transverse measure} $\mu$. We shall refer to the couple $(\FF,\mu)$ as a
\emph{measurable foliation} (see, for instance, \cite{Moore-Schochet88} for
more details).

If a measurable foliation $(\FF,\mu)$ is endowed with a transversely measurable
leafwise Riemannian structure, we shall call it a (leafwise) \emph{Riemannian
measurable foliation}. We shall say that $(\FF,\mu)$ is \emph{hyperbolic} if
$\mu$-a.e.~leaf is $\d$-hyperbolic (with respect to its Riemannian metric). A
particular case is the situation when a.e.~leaf of $\FF$ is a simply connected
manifold with pinched negative sectional curvatures. In this case measurable
systems of boundary measures can be defined as measurable systems of measures
on leafwise unit tangent spheres (by identifying the latter with the leafwise
hyperbolic$\equiv$visibility boundaries via the leafwise exponential maps). In
the general situation one could apply the approach used for hyperbolic
equivalence relations in \secref{sec:continuous bundle} and
\secref{sec:measurable systems}. However, since we are only interested in
\emph{invariant} systems of boundary measures, we shall spare tedious technical
details by reducing foliations to equivalence relations on appropriate systems
of transversals (such a reduction goes back to Plante \cite{Plante75} and Bowen
\cite{Bowen77}).

\begin{defn} \label{def:uniform}
A system $\TT=\{T_i\}$ of transversals of a Riemannian measurable foliation
$\FF$ is \emph{uniform} if
\begin{itemize}
  \item [(i)]
  For any $r>0$ the number of points from $\TT$ in a.e.~leafwise
  ball of radius $r$ is uniformly bounded by a constant $C=C(r)$;
  \item [(ii)]
  The intersection of $\TT$ with a.e.~leaf of $\FF$ is non-empty, and
  there exists a constant $A>0$ such that for any point $x$ of a.e.~leaf its leafwise
  distance to $\FF(x)\cap\TT$ does not exceed $A$.
\end{itemize}
A foliation is \emph{uniform} if it admits a uniform system of transversals.
\end{defn}

Compactness considerations immediately imply

\begin{prop}
Any Riemannian measurable foliation admitting a compact underlying space is
uniform.
\end{prop}

Let now $(\FF,\mu)$ be a uniform Riemannian measurable foliation with a uniform
system of transversals $\TT$. Denote by $R_\FF$ the induced equivalence
relation on $\TT$ (two points are $R_\FF$-equivalent iff they belong to the
same leaf of $\FF$). Then the restriction $\mu_\TT$ of the measure $\mu$ to
$\TT$ is $R_\FF$-quasi-invariant in the sense of \secref{sec:equiv relations}.
Consider on $R_\FF$ the graph structure
$$
K = \{(x,y)\in R_\FF: d(x,y) \le 3A \} \;,
$$
where $d$ denotes the leafwise Riemannian metric, and $A=A(\TT)$ is the
constant from \defref{def:uniform}. Then one can easily see that the leafwise
graph metrics are roughly isometric to the leafwise Riemannian metrics, so
that, in particular, if $\FF$ is hyperbolic, then the graphed equivalence
relation $(\TT,\mu_\TT,R_\FF,K)$ is also hyperbolic, and the leafwise
hyperbolic boundaries with respect to the Riemannian and the graph metrics
coincide (e.g., see \cite[Proposition 7.14]{Ghys-delaHarpe90}).

\begin{prop}
An invariant system of boundary measures on a hyperbolic foliation which is
measurable with respect to one uniform system of transversals is also
measurable with respect to any other uniform system of transversals.
\end{prop}

\begin{proof}
Let $\TT$ and $\ov\TT$ be two uniform systems of transversals, and let $\l$ be
a measurable invariant system of boundary measures of the equivalence relation
$(\TT,R_\FF)$. By \thmref{th:continuous} we may realize $\l$ as the weak$^*$
leafwise limit of measurable systems $\l^n=\{\l_x^n\}$ of interior measures of
$(\TT,R_\FF)$. We shall need the following ``transition maps'' allowing one to
pass from measures on $\TT$ to measures on $\ov\TT$ and back in a measurable
way: denote by $\a_x,\,x\in\TT$ the probability measure uniformly distributed
on the set of points from $\ov\TT$ nearest to $x$ in the leafwise metric of the
leaf $\FF(x)$, and by $\ov\a_{\ov x},\,\ov x\in\ov\TT$ the probability measure
uniformly distributed on the set of points from $\TT$ nearest to $\ov x$. By
the definition of a uniform system of transversals the support of any measure
$\a_x$ lies within a uniformly bounded distance from $x$, and the same property
holds for the measures $\ov\a_{\ov x}$ as well. Obviously, both maps
$x\mapsto\a_x$ and $\ov\x\mapsto\ov\a_{\ov x}$ are measurable. By using the
measures $\a_x$ and $\ov\a_{\ov x}$ one can pass from the sequence $\l^n$ to
the sequence $\ov\l{}^n$ of measurable systems of interior measures of the
equivalence relation $(\ov\TT,R_\FF)$ defined as
$$
\ov\l{}^n_{\ov x}(\ov y)=\sum_{x,y}\ov\a_{\ov x}(x)\,\l_x^n(y)\,\a_y(\ov y) \;.
$$
Since the measures $\l^n_x$ converge in the weak$^*$ topology to a leafwise
constant limit, the measures $\ov\l{}^n_{\ov x}$ also converge to the same
limit, which is therefore measurable with respect to the system of transversals
$\ov\TT$.
\end{proof}

By using the above Proposition, we can now define an invariant measurable
system of boundary measures of a uniform hyperbolic foliation as an invariant
measurable system of boundary measures of any uniform system of transversals.
\thmref{th:at most 2} and \thmref{th:3types} then imply

\begin{thm}
If $\l=\{\l_x\}$ is an invariant measurable system of boundary measures on a
uniform conservative hyperbolic foliation $(\FF,\mu)$, then a.e.~measure $\l_x$
is supported by at most 2 points.
\end{thm}

\begin{thm} \label{thm:classes of foliations}
Let $(\FF,\mu)$ be a uniform ergodic hyperbolic foliation which has more than
one leaf. Then it admits an invariant measurable system of boundary probability
measures if and only if it belongs to one of the following 3 classes:
\begin{itemize}
\item[(i)]
$(\FF,\mu)$ admits a unique invariant measurable system of boundary measures
which is 1-elementary;
\item[(ii)]
$(\FF,\mu)$ admits a unique invariant measurable system of boundary measures
which is 2-elementary;
\item[(iii)]
$(\FF,\mu)$ admits two disjoint 1-elementary invariant measurable systems of
boundary measures, and any invariant measurable system of boundary measures is
their convex combination.
\end{itemize}
\end{thm}

Examples of foliations from each of the classes described in
\thmref{thm:classes of foliations} can be constructed along the same lines as
in \secref{sec:examples}.

\section{Amenability} \label{sec:amenability}

\subsection{Groups}

The class of \emph{amenable groups} is a natural generalization of the class of
finite groups. We shall remind just two definitions of amenability of a locally
compact group $G$ (among a host of other definitions, see \cite{Greenleaf69},
\cite{Paterson88}, \cite{Pier84}). The \emph{Reiter condition}
\begin{equation} \label{eq:reiter}
\ft{\textsf{There exists an {\it approximatively invariant} sequence of
probability measures $\th_n$ on $G$, i.e., such that $\|g\th_n-\th_n\|\to 0$
for any $g\in G$.}}
\end{equation}
by its constructiveness is very useful for verifying amenability of a given
group, whereas, on the contrary, the \emph{fixed point condition}
\begin{equation} \label{eq:fixed}
\ft{\textsf{Any continuous affine action of $G$ on a compact space has a fixed
point.}}
\end{equation}
is the one which is most useful for applications of amenability. For instance,
the fixed point condition implies at once that any continuous action of an
amenable group on a compact space has an invariant measure.

We shall illustrate how these two definitions work together on the example of
the following result. Recall that a group of isometries of a hyperbolic space
$\XX$ is called \emph{elementary} if it either is compact or fixes a finite
subset of $\part\XX$ (actually, if in the latter case the fixed set contains
more than 2 points, then the group fixes also a compact subset of $\XX$, and
therefore is compact, see \propref{pr:allbary}).

\begin{thm} \label{th:amgroup}
A closed subgroup $G$ of the group of isometries of an exponentially bounded
hyperbolic space $\XX$ is amenable if and only if it is elementary.
\end{thm}

\begin{proof}
If $G$ is amenable, then by the fixed point property \eqref{eq:fixed} it has an
invariant probability measure $\l$ on $\part\XX$. If $\l$ does not have atoms
of weight at least $\displaystyle\frac12$, then by \propref{pr:allbary} its
quasi-barycenter set is compact and fixed by $G$ together with the measure
$\l$. Otherwise, $G$ fixes the set of maximal weight atoms of the measure $\l$
which consists either of 1 or of 2 points.

Conversely, if $G$ is compact, then it is obviously amenable. Otherwise, if $G$
fixes a point $\g\in\part\XX$, take a point $o\in\XX$ and put
$\th_n=\l_n(o,\g)$, where $\l_n$ are the maps from \thmref{th:lambda}. Then
$$
\|g\th_n-\th_n\| = \|g\l_n(o,\g)-\l_n(o,\g)\| = \|\l_n(go,\g)-\l_n(o,\g)\| \to
0 \;,
$$
so that Reiter's condition \eqref{eq:reiter} is satisfied. If the fixed set of
$G$ consists of two points $\{\g_1,\g_2\}$, then just take
$\th_n=\bigl[\l_n(o,\g_1)+\l_n(o,\g_2)\bigr]/2$.
\end{proof}

This result is classical for the constant curvature spaces (i.e., for Kleinian
groups), and, of course, should be of no surprise in the stated generality
either, although I could not find it in the literature in an explicit form. For
the case of a group with a finite critical exponent acting on a CAT($-1$) space
it is proved in \cite[Corollary 2.5]{Burger-Mozes96} (also see
\cite{Burger-Schroeder87} and \cite{Adams-Ballmann98} for properties of
amenable groups of isometries of CAT($0$) spaces); note that the assumption of
exponential boundedness of the space $\XX$ in \thmref{th:amgroup} is essential
(see \secref{sec:counterexamples} for the corresponding counterexamples).

However, the main reason for giving this proof is that it contains the key
geometric ingredients of the relationship between amenability and
hyperbolicity, namely, boundary convergence of geodesics (essentially
equivalent to \thmref{th:lambda}) and existence of quasi-barycenter sets, as
illustrated by the following diagram:

\bigskip

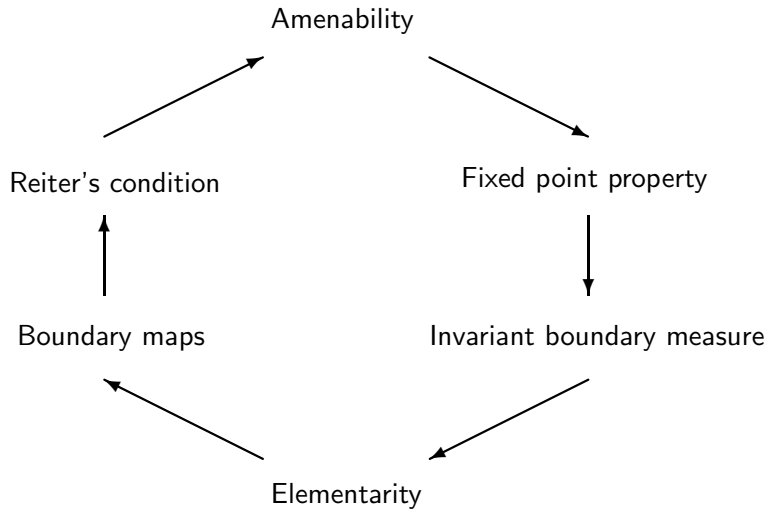
\setlength{\unitlength}{1pt}
\begin{picture}(340,200)(0,0)
\thicklines

\put(120,180){\makebox(0,0)[lb]{\large\textsf{Amenability}}}

\put(180,171){\vector(2,-1){60}}

\put(57,141){\vector(2,1){60}}

\put(21,120){\makebox(0,0)[lb]{\large\textsf{Reiter's condition}}}

\put(192,120){\makebox(0,0)[lb]{\large\textsf{Fixed point property}}}

\put(240,111){\vector(0,-1){30}}

\put(57,81){\vector(0,1){30}}

\put(24,60){\makebox(0,0)[lb]{\large\textsf{Boundary maps}}}

\put(180,60){\makebox(0,0)[lb]{\large\textsf{Invariant boundary measure}}}

\put(117,19){\vector(-2,1){60}}

\put(240,49){\vector(-2,-1){60}}

\put(120,0){\makebox(0,0)[lb]{\large\textsf{Elementarity}}}

\end{picture}

\bigskip

The notion of amenability can also be defined for a number of objects other
than groups: discrete equivalence relations, foliations, groups actions, see
\cite{Zimmer78}, \cite{Connes-Feldman-Weiss81}, \cite{Hurder-Katok87},
\cite{Adams92}, \cite{Adams-Elliott-Giordano94}, \cite{Kaimanovich97} for
various definitions. All these definitions can be put into the general
framework of amenability for \emph{measured groupoids}, see
\cite{Anantharaman-Renault00}, \cite{Corlette-Lamoneda-Iozzi02}. We shall refer
the reader to the aforementioned references for precise definitions; note that
in all these cases it is possible to give (equivalent) definitions of
amenability both in terms of approximatively invariant sequences of probability
measures (analogous to the Reiter condition \eqref{eq:reiter}; we shall briefly
mention these definitions when dealing with the corresponding classes of
objects) and in terms of an appropriate fixed point property (analogous to
condition \eqref{eq:fixed}).

We shall now proceed to establish analogues of \thmref{th:amgroup} for these
objects by taking stock of the results from \secref{sec:hyperbolic} and
\secref{sec:relations}, and following precisely the same scheme: amenability by
the fixed point property implies existence of an invariant system of boundary
measures, and therefore elementarity; whereas elementarity by virtue of
\thmref{th:lambda} leads to existence of an approximatively invariant sequence
of probability measures, i.e., to amenability.

\subsection{Equivalence relations}

\begin{thm} \label{th:equivmain}
A conservative hyperbolic equivalence relation with uniformly bounded vertex
degrees is amenable if and only if it is elementary.
\end{thm}

\begin{proof}
If the equivalence relation is amenable, then the fixed point property applied
to the Banach bundle of leafwise boundary measures implies existence of an
invariant measurable system of boundary measures, see \cite[Theorem
4.2.7]{Anantharaman-Renault00} (based on the original argument from
\cite{Zimmer77}), which by \thmref{th:at most 2} implies that the equivalence
relation is elementary.

For proving the converse implication we have to construct a sequence of
measurable maps $\l_n$ assigning to a.e. point $x$ from the state space of the
equivalence relation a probability measure on the equivalence class $[x]$, and
such that for a.e. pair of equivalent points $(x,y)$
\begin{equation} \label{eq:convergence}
\| \l_n(x) - \l_n(y) \| \to 0
\end{equation}
(this is the analogue of the Reiter condition for equivalence relations).

The only point which needs an explanation here is the measurability of the
approximatively invariant sequences of measures obtained in \thmref{th:lambda}.
It would follow from the existence of a measurable family of leafwise tempered
measures. Under the assumption of uniform boundedness of vertex degrees the
leafwise counting measures are obviously tempered.
\end{proof}

\begin{rem}
S. Adams \cite{Adams90} proved that a treed equivalence relation with a
\emph{finite invariant measure} is amenable if and only if a.e.~tree has at
most 2 ends. It is plausible that the same is true for hyperbolic equivalence
relations as well, namely, that for a hyperbolic equivalence relation with a
finite invariant measure the hyperbolic boundary of a.e.~leaf has at most 2
points.
\end{rem}

\begin{rem}
It is well-known that if a graphed equivalence relation has
\emph{subexponential growth} (i.e., $[\card B(x,n)]^{1/n}\to 1$ for a.e. $x\in
X$, where $B(x,n)$ are the balls in the leafwise graph metrics), then it is
amenable, e.g., see \cite[Proposition~3.3]{Adams-Lyons91}. Since the sequences
of leafwise balls centered at any two equivalent points are obviously
sandwiched, \lemref{lem:sandwich} provides a direct ``constructive'' way of
proving this fact, and this Lemma can also be used in numerous other arguments
deducing amenability from subexponentiality in various contexts, see
\cite{Adelson-Shreider57}, \cite{Greenleaf69}, \cite{Series79},
\cite{Samuelides79}, \cite{Anantharaman-Renault00},
\cite{delaHarpe-Grigrochuk-Ceccherini99}. However, in \thmref{th:lambda}
\lemref{lem:sandwich} is used in a somewhat different situation, where the
sandwiched sequences of sets are not metric balls.
\end{rem}

\begin{rem}
The boundary maps provided by \thmref{th:lambda} do not depend on the measure
on the state space of the equivalence relation, so that in fact the argument in
the proof of \thmref{th:equivmain} shows amenability of elementary hyperbolic
equivalence relations even in the \emph{Borel} (rather than measure) category
(cf. the discussion of definitions of amenable action in various categories in
\secref{sec:actions}).
\end{rem}

\subsection{Foliations}

In view of the discussion in \secref{sec:fol}, measurable foliations can be
dealt with in precisely the same way as equivalence relations (by using
\thmref{th:lambda} and \thmref{thm:classes of foliations}). Once again, one of
equivalent definitions of amenability amounts to existence of an
approximatively invariant sequences of measures \eqref{eq:convergence} (note
that there is also a somewhat different notion of amenable foliation meaning
that its \emph{fundamental groupoid} is amenable, see the instructive
discussion in \cite{Corlette-Lamoneda-Iozzi02}).

\begin{thm} \label{th:foliations main}
A uniform conservative hyperbolic foliation is amenable if and only if it is
elementary.
\end{thm}

\begin{cor}
A hyperbolic foliation admitting a compact underlying space is amenable if and
only if it is elementary.
\end{cor}

\begin{ex}
Let $\FF$ be a foliation with simply connected leaves of pinched negative
sectional curvature. Denote by $S\FF$ the associated \emph{stable extension}.
This is the foliation whose underlying space consists of leafwise unit tangent
vectors of $\FF$, and the leaves consist of tangent vectors from the same leaf
of $\FF$ pointing at the same point at infinity. Then the foliation $S\FF$ is
obviously elementary, so that it is amenable with respect to any
quasi-invariant transverse measure. This result is a generalization of the
well-known amenability of the stable foliation of the geodesic flow, see
\cite{Bowen77}. Yet another example of elementary hyperbolic foliations is
provided by foliations and laminations with pointed at infinity leaves arising
in conformal dynamics \cite{Kaimanovich-Lyubich01}.
\end{ex}

\begin{rem}
The construction of the stable extension can also be applied to an arbitrary
hyperbolic foliation (not necessarily with Cartan--Hadamard leaves), in which
case we shall call it the \emph{boundary extension}. The underlying space of
the boundary extension consists of the pairs $(x,\g)$, where $x$ is a point of
the underlying space of the original foliation, and $\g$ is a point of the
hyperbolic boundary of the leaf of $x$. The leaves of the boundary extension
consist of all pairs $(x,\g)$ with the same $\g$. The boundary extension is
elementary and therefore amenable.
\end{rem}

\subsection{Group actions} \label{sec:actions}

Let $G$ be a locally compact group with a measure type preserving action on a
measure space $(X,\mu)$. The ``Reiter type'' definition of amenability of this
action requires existence of a sequence of measurable maps $\l_n:X\to\Pr(G)$
which are \emph{approximatively equivariant} in the sense that
\begin{equation} \label{eq:amenable action}
\| \l_n(gx) - g\l_n(x) \| \to 0
\end{equation}
for all $g\in G$ and $\mu$-a.e. $x\in X$. There are several modifications of
this definition. In particular, a Borel action of a locally compact group $G$
on a Borel space $X$ is called \emph{universally amenable} \cite{Adams96}
$($other terms: \emph{measurewise amenable} \cite{Anantharaman-Renault00},
\emph{measure-amenable} \cite{Jackson-Kechris-Louveau02}$)$ if it is amenable
for any quasi-invariant measure $\mu$ on $X$. If, moreover, $X$ is a
topological space, and the action of $G$ on $X$ is continuous, then the action
is called \emph{topologically amenable} if the convergence in
\eqref{eq:amenable action} is uniform on compact subsets of $G\times X$.

\begin{thm} \label{th:actions mains}
An ergodic action of a word hyperbolic group $G$ on a measure space $(X,\mu)$
is amenable if and only if it factorizes through the action of $G$ either on
the hyperbolic boundary $\part G$, or on its symmetric square $\wt{\part^2 G}$.
\end{thm}

\begin{proof}
The amenability of the action of $G$ on $\part G$ or on $\wt{\part^2 G}$ is
well-known (see \remref{rem:bdry am} below) and follows at once from
\thmref{th:lambda}. Therefore, any action which factorizes through $\part G$ or
$\wt{\part^2 G}$ is also amenable.

Conversely, amenability of the action on the space $(X,\mu)$ by the fixed point
property applied to the Banach bundle of boundary measures implies existence of
an equivariant measurable map from $X$ to the space of probability measures on
$\part G$ (in precisely the same way as in the proof of \thmref{th:equivmain}).
By ergodicity we may assume that either these measures are a.e. purely
non-atomic, or that almost all of them are equidistributed on a set of atoms
whose cardinality is a.e. the same. In the first case, or if the number of
atoms exceeds 2, we obtain by \propref{pr:allbary} an equivariant map from $X$
to the space of finite subsets of $G$, which means that the action on $X$ is
dissipative, and therefore is confined to a single orbit. Otherwise we have the
sought for equivariant projection from $X$ to $\part G$ or $\wt{\part^2 G}$.
\end{proof}

\begin{thm} \label{th:univ}
Let $G$ be a closed subgroup of the group of isometries of an exponentially
bounded hyperbolic space $\XX$. Then the action of $G$ on the hyperbolic
boundary $\part\XX$ is topologically amenable.
\end{thm}

\begin{proof}
A sequence of approximatively invariant measures for this action is provided by
\thmref{th:lambda}. Indeed, let $\l_n:\XX\times\part \XX\to\Pr(\XX)$ be the
$\Iso(\XX)$-equivariant maps constructed in \thmref{th:lambda}. Choose a
reference point $o\in\XX$, and put $\wt\l_n(\g)=\l_n(o,\g)$ for any $g\in G$
and $\g\in\part\XX$. Then the maps $\wt\l_n:\part\XX\to\Pr(\XX)$ are
approximatively equivariant, because
$$
\| \wt\l_n(g\g) - g\wt\l_n(\g) \| = \| \l_n(o,g\g) - \l_n(go,g\g) \| \to 0
$$
for any $g\in G$ and $\g\in\part\XX$ uniformly on compacts, see
\remref{rem:uniform}. Since the action of $G$ on $\XX$ is proper, there exists
an equivariant Borel lift $\Pr(\XX)\to\Pr(G)$ (e.g., see
\cite[Proposition~2.1.10]{Anantharaman-Renault00}). Applying it to the sequence
$\wt\l_n$ gives an approximatively equivariant sequence of maps from $\part\XX$
to $\Pr(G)$.
\end{proof}

\begin{rem} \label{rem:bdry am}
The action of the group $SL(2,\Z)$ on the circle at infinity of the hyperbolic
plane was the first example of an amenable action of a non-amenable group
\cite{Bowen77}, \cite{Vershik78}. Since then amenability of boundary actions
was established in a number of other situations \cite{Spatzier87},
\cite{Spatzier-Zimmer91}, \cite{Adams94}; the most general result being
universal amenability of the boundary action of the group of isometries of an
exponentially bounded hyperbolic space \cite{Adams96}. In the particular case
of the boundary action of a word hyperbolic group the arguments of Adams were
recently simplified by Germain \cite{Germain00}, who proved topological
amenability of this action. However, the approach based on using
\thmref{th:lambda} and \lemref{lem:sandwich} is more streamlined even in this
case.
\end{rem}

In the CAT($-1$) case replacing \thmref{th:lambda} with \thmref{th:lambda-CAT}
in the proof of \thmref{th:univ} yields

\begin{thm} \label{th:univ-CAT}
Let $G$ be a closed subgroup of the group of isometries of a hyperbolic space
$\XX$ with finite critical exponent. Then the action of $G$ on the hyperbolic
boundary $\part\XX$ is topologically amenable.
\end{thm}

\begin{rem}
Burger and Mozes in \cite[Corollary 1.4]{Burger-Mozes96} stopped short of
saying that under the assumptions of \thmref{th:univ-CAT} the action of $G$ on
$\part\XX$ is amenable with respect to any Patterson measure (although they
basically proved that). Moreover, they also made the first step to establishing
universal amenability of the boundary action by showing that the stabilizer of
\emph{any} boundary point is amenable (which is a necessary condition for
universal amenability of the action, see \cite{Adams92}).
\end{rem}

\begin{rem}
As we shall see in \secref{sec:counterexamples} below, without any assumptions
on the growth of the space $\XX$ \thmref{th:univ} and \thmref{th:univ-CAT} may
fail even in the weakest form of amenability of the action with respect to a
given quasi-invariant boundary measure, which answers a question formulated by
Adams in \cite{Adams96}.
\end{rem}

\subsection{Counterexamples to boundary amenability}
\label{sec:counterexamples}

As it was pointed out by Burger and Mozes in \cite{Burger-Mozes96} (Remark
after Proposition 1.6; in our notations $X(\infty)\equiv\part X$),
\begin{quotation}\small
The product $X = [0,\infty)\times\H^2$ with Riemannian metric $(dt)^2 +
e^{-2t}\frac{|dz|^2}{y^2}$, where $\H^2$ denotes the upper half-plane, is a
CAT$(-1)$ space on which $SL(2,\R)$ acts isometrically, fixing a point in
$X(\infty)$.
\end{quotation}
Note that the metric on $X$ is obtained from the hyperbolic metric on $\H^2$ in
precisely the same way as the hyperbolic metric on (a horoball in) $\H^2$ is
obtained from the Euclidean metric on $\R$, i.e., by exponential rescaling in
the vertical direction. In fact, this example is a particular case of the
following general \emph{hyperbolization construction} (cf.
\cite{Kaimanovich-Lyubich01}). This construction resembles that of a
$\kappa$-\emph{cone over a metric space} (e.g., see
\cite{Bridson-Haefliger99}).

\begin{prop} \label{pr:extension}
For a complete \emph{CAT$(0)$} space $Y$ put $X=\R\times Y$, and define the
metric $d_X$ on $X$ by putting for any two points $(t_1,y_1),(t_2,y_2)\in X$
\begin{equation} \label{eq:dX}
d_X \bigl( (t_1,y_1), (t_2,y_2) \bigr) = d_{\H^2} \bigl(
(t_1,0),(t_2,d_Y(y_1,y_2)) \bigr) \;,
\end{equation}
where $d_Y$ is the metric on $Y$, and $d_{\H^2}$ is the metric
$$
ds^2 = dt^2 + e^{-2t} dl^2
$$
on the hyperbolic plane $\H^2\cong\R\times\R$ realized in the upper half-plane
model with logarithmic rescaling of the vertical coordinate. Then $d_X$ is a
complete metric on $X$, the space $X$ is \emph{CAT$(-1)$} with respect to this
metric, and the boundary $\part X$ of the hyperbolic compactification of $X$ is
the union $Y\cup\{\om\}$, where $\om$ is the common limit point of vertical
geodesic rays $[0,\infty)\times\{y\},\,y\in Y$ in $X$ directed upwards. The
point of $\part X$ corresponding to a point $y\in Y$ is the limit point of the
directed downwards vertical geodesic ray $[0,-\infty)\times\{y\}$. If $Y$ is
non-compact, then $\part X$ is homeomorphic to the one-point compactification
of $Y$, and if $Y$ is compact, then $\part Y$ is homeomorphic to the disjoint
union of $Y$ and $\om$.
\end{prop}

\begin{defn}
We shall call the CAT$(-1)$ space $(X,d_X)=(\HH Y,d_{\HH Y})$ constructed in
\propref{pr:extension} the \emph{hyperbolization} of the CAT$(0)$ space
$(Y,d_Y)$.
\end{defn}

\begin{proof}[Proof of \propref{pr:extension}]
The verification of the triangle inequality for $d_X$ using formula
\eqref{eq:dX} is straightforward. Completeness of the metric $d_X$ is obvious.
Further, the space $(X,d_X)$ is geodesic, and the (unique) geodesic $\ov\xi$
joining two points $(t_1,y_1)$ and $(t_2,y_2)$ is the image of the geodesic
joining the points $(t_1,0)$ and $(t_2,d_Y(y_1,y_2))$ in $\H^2$ under the
natural isometry between the strip $\R\times[0,d_y(y_1,y_2)]$ in $\H^2$ and the
``curtain'' $\R\times\xi$ in $X$ over the geodesic $\xi$ joining $y_1$ and
$y_2$ in $Y$.

In order to check the CAT$(-1)$ property let us consider a geodesic triangle
$\D_X=\bigl[ (t_0,y_0),(t_1,y_1),(t_2,y_2) \bigr]$ in $X$. It is more
convenient to construct a comparison triangle for $\D_X$ in the 3-dimensional
hyperbolic space $\H^3$ rather than in the 2-dimensional hyperbolic plane
$\H^2$, and to use the fact that $\H^3$ is obtained from the Euclidean plane
$\E^2$ in precisely the same way as $X$ is obtained from $Y$. Let
$\D_{\E^2}=[z_0,z_1,z_2]$ be the comparison triangle in $\E^2$ for the triangle
$\D_Y=[y_0,y_1,y_2]$ in $Y$. Then the definition of the distance $d_X$ implies
that $\D_{\H^3}=\bigl[ (t_0,z_0),(t_1,z_1),(t_2,z_2) \bigr]$ is a comparison
triangle for $\D_X$. Now take points $(t'_i,y'_i),\,i=1,2,$ on the sides
$[(t_0,y_0),(t_i,y_i)]$ of the triangle $\D_X$, and denote by $(t'_i,z'_i)$ the
corresponding points on the sides of the comparison triangle $\D_{\H^3}$. By
the above description of the geodesics in $X$ the ``heights'' $t'_i$ remain the
same, and $d_Y(y_0,y'_i)=d_{\E^2}(z_0,z'_i)$. Therefore the CAT$(0)$ property
of the space $Y$ and the definition of the metric $d_X$ imply that
$$
d_X\bigl( (t'_1,y'_1), (t'_2,y'_2) \bigr)
 \le d_{\H^3}\bigl( (t'_1,z'_1),(t'_2,z'_2) \bigr) \;.
$$

The description of the hyperbolic boundary $\part X$ easily follows from the
above description of the geodesics in $X$.
\end{proof}

\propref{pr:extension} now allows one to obtain examples showing that
\thmref{th:amgroup}, \thmref{th:univ} and \thmref{th:univ-CAT} may well fail
without any bounded geometry assumptions on the space or on the group.

\begin{prop} \label{pr:counter}
Let $Y$ be any \emph{CAT$(0)$} space with a non-amenable closed group of
isometries $G$, and let $X=\HH Y$ be its hyperbolization as described in
\propref{pr:extension}. Then the group $G$ acting on $X$ as $g(t,y)=(t,gy)$ is
also closed in $\Iso(X)$, and it fixes the point $\om\in\part X$. Therefore,
\begin{itemize}
\item
The group $G$ is elementary with respect to the space $X$ and non-amenable;
\item
The action of the group $G$ on $\part X$ is non-amenable with respect to the
invariant measure $\d_\om$ concentrated at the fixed point $\om$.
\end{itemize}
\end{prop}

In order to obtain more concrete examples one can take $Y$ to be, for instance,
the hyperbolic plane (as did Burger and Mozes; actually they took just the
horoball in $\HH\H^2$ centered at the point $\om$ rather than the whole space
$\HH\H^2$), or the Cayley 1-complex of the free group. Still, it would be
interesting to have examples of non-amenable actions with respect to a
non-atomic boundary measure.

\bibliographystyle{amsalpha}

\providecommand{\bysame}{\leavevmode\hbox to3em{\hrulefill}\thinspace}
\providecommand{\MR}{\relax\ifhmode\unskip\space\fi MR }
\providecommand{\MRhref}[2]{%
  \href{http://www.ams.org/mathscinet-getitem?mr=#1}{#2}
} \providecommand{\href}[2]{#2}

\end{document}